\documentclass[11pt]{amsart}
\usepackage{amssymb}
\usepackage{multirow}

\usepackage{palatino}

\newtheorem{thm}{Theorem}[section]
\newtheorem{cor}[thm]{Corollary}
\newtheorem{prop}[thm]{Proposition}
\newtheorem{lem}[thm]{Lemma}

\theoremstyle{remark}
\newtheorem{rmk}[thm]{Remark}

\newtheorem{quest}[thm]{Question}

\newcommand{\Lg}{\mbox{$\mathfrak{g}$}}
\newcommand{\Lk}{\mbox{$\mathfrak{k}$}}

\newcommand{\Lh}{\mbox{$\mathfrak{h}$}}

\newcommand{\Ld}{\mbox{$\mathfrak{d}$}}

\newcommand{\Ja}{\mbox{$\mathfrak{A}$}}
\newcommand{\Js}{\mbox{$\mathfrak{S}$}}
\newcommand{\Jn}{\mbox{$\mathfrak{N}$}}
\newcommand{\Jb}{\mbox{$\mathfrak{B}$}}
\newcommand{\Jt}{\mbox{$\mathfrak{T}$}}

\newcommand{\Ad}{\mbox{Ad}}

\newcommand{\Pf}{{\em Proof}. }
\newcommand{\EPf}{\hfill$\square$}

\newcommand{\R}{\mbox{$\mathbb R$}}
\newcommand{\C}{\mbox{$\mathbb C$}}

\newcommand{\bigzero}{\mbox{\normalfont\Large\bfseries 0}}


\begin{document}


\title{A moment map for the variety of Jordan algebras}


\author{Claudio Gorodski and Iryna Kashuba and Mar\'\i a Eugenia Martin}
\thanks{The first author acknowledges partial financial
  support from  CNPq (grant 304252/2021-2) and FAPESP (grant 16/23746-6).
The second author acknowledges partial financial support from FAPESP (grant 2021/14716-4).}

\date{\today} 

\address{}
\email{}

\begin{abstract}
  We study the variety of complex $n$-dimensional Jordan algebras
  using techniques from Geometric Invariant Theory. More specifically,
  we use the Kirwan-Ness theorem to construct a Morse-type stratification 
  of the variety of Jordan algebras into finitely many invariant smooth
  subvarieties, with respect to the energy functional associated to
  the canonical moment map. In particular we obtain a new, cohomology-free
  proof of the well-known rigidity of semisimple Jordan algebras. 
\end{abstract}

\maketitle 

\section{Introduction}\label{intro}

The main goal of this paper is to introduce techniques from
Geometric Invariant Theory (GIT) in the study of finite-dimensional
complex Jordan algebras. In particular, we consider a naturally
defined moment map for the affine variety of $n$-dimensional
Jordan algebras and the associated energy functional, and  
use it to define a notion of soliton Jordan algebra.
Somehow in a dual sense, the search for a soliton
Jordan algebra in an isomorphism class of Jordan
algebras can also be thought of as the search for a 'best' (Hermitean)
metric in a Jordan algebra. Herein we also make an effort to relate the
new geometric invariants to the 'old' algebraic invariants of Jordan algebras
(with partial success).

\subsection{Basic setup}

In this work we only consider finite-dimensional algebras over $\C$.  
A Jordan algebra is a commutative algebra $\Ja$ satisfying the following
Jordan identity:
\[ [L_a,L_{a^2}]=0 \]
for all $a\in\Ja$, where $L_a:\Ja\to\Ja$ denotes the left-multiplication 
by~$a$, and $[,]$ denotes the commutator. Here we only note that
Jordan algebras are, in general, non-associative, and refer the reader
to~\cite{mccrimmon} for more about Jordan algebras and their 
importance.

We can view a commutative algebra of dimension~$n$ as a 
commutative multiplication in $\C^n$, namely, as a tensor
$\mu\in S^2(\C^n)^*\otimes\C^n\cong S^2(\C^{n*})\otimes\C^n=:V_n$
(symmetric in the first two arguments). 
Further, the Jordan identity is equivalent to its linearized form
\[ (ab,c,d)+(bd,c,a)+(da,c,b) = 0, \]
for all $a$, $b$, $c$, $d\in\Ja$, where 
$(a,b,c)=(ab)c-a(bc)$ is the associator of $a$, $b$, $c$.
These are polynomial equations in $V_n$, therefore the space $\mathcal J_n$
of Jordan algebras of a dimension $n$
can be naturally identified with an affine variety in $V_n$.
It is useful to note that these polynomials are homogeneous, so we could
also view the space of Jordan algebras modulo complex scaling
as a projective variety in $\mathbb PV_n$; however, we prefer to do
calculations in the vector space $V_n$. 

There is a natural action of the group $G:=GL(n,\C)$ on the space
$V_n$, which corresponds to 'change of basis'. Namely, for 
$g\in G$ and $\mu\in V_n$, one puts
\begin{equation}\label{g-action}
  g\cdot \mu (a,b)=g(\mu(g^{-1}a,g^{-1}b))
\end{equation}
for all $a$, $b\in \C^n$. It is clear that the orbit $G\cdot\mu$ yields
the isomorphism class of $\mu$. In particular, $\mathcal J_n$ is an
invariant subvariety.

In order to use methods from GIT, we need to fix a background
Hermitean product $Q$ on $\C^n$
(whereas in the case of Lie algebras
$Q$ would correspond to a left-invariant Hermitean metric on the
corresponding simply-connected complex Lie group, in our case there
is no obvious global object attached to $Q$). However, we will see
that the results that we shall obtain for isomorphism
classes of Jordan algebra will not depend on $Q$ in an essential way.
The group of unitary transformations of $(\C^n,Q)$ is a
maximal compact subgroup $K\cong U(n)$ of $G$, which also
acts by unitary transformations on $V_n$, if this space
is equipped with the induced Hermitean product
(which we denote with the same symbol).
Note that fixing~$Q$ is equivalent to fixing~$K$.  
Now, it is standard that the
action of the \emph{compact} group $K$ on the projective variety
$\mathbb PV_n$ yields a moment map $\mathbb PV_n\to \Lk^*$,
where $\Lk$ is the Lie algebra of $K$ and the star refers to the dual
space; we prefer to do calculations in the vector space $V_n$
and we use a different normalization, so
we shall equivalently consider a 
scale invariant moment map $m: V_n\setminus\{0\}\to i\Lk$
(see below for precise definitions);
note that $i\Lk$ is the space of Hermitean transformations of $(\C^n,Q)$.
The square norm $E_n=||m||^2$, which we call \emph{energy}, has
the very interesting property that its gradient flow lines 
are contained in the orbits of $G$, so the limits of those lines
define degenerations of the initial algebra.
This brings us to the study of critical
points of $E_n$, which we call \emph{solitons} (for reasons that
shall become clear later). 

The Kempf-Ness theorem~\cite{kempf-ness}
provides a beautiful characterization of
closed orbits of a rational representation $V$ of a complex reductive
linear algebraic group $G$, namely, these are precisely
the orbits that contain a
\emph{minimal} vector, i.~e.~a vector realizing the distance of the
orbit to the origin (with respect to a fixed Hermitean product~$Q$).  
Since the moment map associated to~$Q$
controls the infinitesimal change of the norm of vectors,
it turns out that minimal vectors precisely correspond
to zeros of the moment map (or the energy functional). 
The Kempf-Ness theorem also implies that a non-closed
orbit with positive distance to the
origin contains a non-trivial closed orbit in its closure.
In particular, this gives a description of the categorical
quotient $V/\!/G$ (which is by definition the affine variety defined by the
$G$-invariant polynomials in $V$). 

On the other hand, the union of the orbits containing
zero in its closure is called the \emph{null cone} $N$
of the representation, 
and its elements are called \emph{unstable} vectors.
In our case of $GL(n,\C)$ acting on $V_n$, every vector
is unstable since we can always rescale the multiplication that
defines the algebra, which is to say that $N=V_n$.
Nonetheless, Kirwan~\cite{kirwan} and Ness~\cite{ness}
showed how to use the moment map to study the orbit space
of $N$, for a general representation $V$ of a
complex reductive linear algebraic group $G$;
herein we are not so much interested in this
orbit space, but in the related stratification of $N$.
Consider a Hermitian product on $V$, the associated
maximal compact subgroup $K$ of $G$, the
associated moment map, energy, and solitons.
Kirwan and Ness proved that all non-minimal
solitons, i.~e.~those solitons with positive energy, do occur in $N$.
Further, there is only a finite set of distinguished $G$-orbits in $N$
which contain solitons; moreover, for each such orbit $G\cdot v$,
the set of solitons comprise a single $K$-orbit (up to scaling)
and they are minima of the energy on~$G\cdot v$. Finally, 
although the energy is not a Morse(-Bott) function,
the Kirwan-Ness theorem yields a Morse-type
stratification of $N$ into finitely many
invariant smooth subvarieties, each of which being the stable
manifold of the set of solitons of a certain \emph{type}. 

In this paper, we adapt such GIT methods to  
the setting of Jordan algebras, in parallel with some of  
the work  of Lauret for Lie algebras~\cite{lauret2003}, and use them to 
prove a number results which we explain in the sequel. We also
formulate some conjectures and open problems.  

In the following, in case there is no danger of confusion,
we denote $Q$ simply by
$\langle\cdot,\cdot\rangle$ and the associated norm by $||\cdot||$.

\subsection{Main results}
Recall that a \emph{(geometric) deformation} of a Jordan algebra 
$\mu\in\mathcal J_n$ is a Jordan algebra $\nu\in\mathcal J_n$ 
whose $G$-orbit contains $\mu$ in its Zariski closure
(which is the same as closure in the Hausdorff topology),
that is $\mu\in \overline{G\cdot\nu}$. 
In this case, we also say that $\mu$ is a 
\emph{degeneration} of $\nu$ and we write $\nu\to\mu$.
The degeneration $\nu\to\mu$ is called \emph{trivial} if
$\mu$ is isomorphic to $\nu$. 
A Jordan algebra $\mu\in\mathcal J_n$
is called \emph{(geometrically) rigid} if its
$G$-orbit is Zariski-open; in our context (finite-dimension,
base field $\C$), this is equivalent to the non-existence of 
non-trivial deformations of $\mu$ or, yet, to having the closure of
the orbit of $\mu$ in $\mathcal J_n$ coincide with an irreducible
component of $\mathcal J_n$.  

A commutative algebra $\mu\in V_n$ will be called a \emph{soliton}
if $\mu$ is a critical point of $E_n$. 
A $G$-orbit in~$V_n$ (or isomorphism class of
commutative algebras) is called \emph{distinguished}
if it contains a soliton.
For $\mu\in V_n\setminus\{0\}$, we write 
$M_\mu=||\mu||^2m(\mu)$ and call it the \emph{moment matrix} of $\mu$.
We have a structure result for Jordan solitons.

\begin{prop}\label{prop:struct}
Let $\mu\in V_n$ and write $\Ja=(\C^n,\mu)$ 
for the corresponding commutative algebra. Then:
\begin{enumerate}
\item[(a)] $\mu$ is a soliton if and only if 
its moment matrix
$M_\mu=c_\mu I+D_\mu$, where $c_\mu<0$, $I$ denotes the identity matrix,
and $D_\mu$ is a derivation of $\mu$ and a Hermitean matrix with 
respect to~$Q$; further, in this case there is a positive multiple
of $D_\mu$ which has rational eigenvalues.  
\item[(b)] If $\mu$ is a Jordan soliton 
then there is a maximal semisimple subalgebra $\Js$ of $\Ja$ 
such that $\Ja=\Js+\Jn$ is a $D_\mu$-invariant decomposition, where 
$\Jn$ denotes the radical of $\Ja$.  Further,
$D_\mu|\Js=0$. 
\end{enumerate}
\end{prop}

Using this structure result
and the classic theory of semisimple
Jordan algebras (see~\cite{braun-koecher}), we prove:

\begin{thm}\label{ss-distinguished}
  Every isomorphism class of
  complex semisimple Jordan algebras is distinguished. 
Further, in each dimension the semisimple Jordan solitons 
are characterized among Jordan algebras as having the 
lowest possible value of the energy for that dimension.    
\end{thm}

The relation $\mu\to\nu$ defines a partial order on $\mathcal J_n$, and 
it follows from properties of the stratification that
the functional $E_n$ behaves well with regard to this partial order:

\begin{prop}\label{deg-decr-e}
  Let $\mu$, $\nu\in\mathcal J_n$. If $G\cdot\mu$ is a distinguished
  orbit and the energy of a soliton in $G\cdot\mu$ is larger than
the energy of $\nu$, then 
$\mu$ cannot degenerate to $\nu$. 
\end{prop}

We apply Theorem~\ref{ss-distinguished}
to give a simple, cohomology-free proof of the 
following well-known result~\cite{glassman,finston}. 

\begin{thm}[Glassman, Finston]\label{ss}
Every finite-dimensional complex semisimple Jordan algebra is rigid. 
\end{thm}

The class of nilpotent Jordan algebras of dimension $n$ forms a 
$G$-invariant closed subset $\mathcal N_n$ of~$\mathcal J_n$. 
It is  known that there is no
non-trivial $n$-dimensional complex
Jordan algebra lying in all irreducible components of $\mathcal J_n$
(see~\cite[Theorem~4.65]{martin}). 
On the other hand, by studying some specific degenerations, we 
prove: 

\begin{thm}\label{nilp-comp}
For all $n\geq2$ there exists a non-trivial $n$-dimensional complex
nilpotent Jordan algebra lying in all irreducible components 
of $\mathcal N_n$; further, the orbit of this algebra is characterized
as realizing the highest possible value of the energy~$E_n$.
It follows that the projectivization $\mathbb P\mathcal N_n$ is connected
for $n\geq2$. 
\end{thm}

Finally, we give a complete explicit description of the case $n\leq4$,
see the tables below.
The \emph{type} of a soliton $\mu$ is $(d_1<\cdots<d_r;m_1,\ldots,m_r)$,
where $d_1,\ldots,d_r$ is the sequence of coprime integers which is
proportional to the eigenvalues of $m(\mu)+E_n(\mu)I=\frac1{||\mu||^2}D_\mu$,
and $m_1,\ldots,m_r$ are the respective multiplicities.  
The stratification is finite and the strata are parametrized
by~$\beta=m(\mu)$ or, equivalently, by the soliton type.
We deduce:

\begin{thm}\label{low}
  Every isomorphism class of complex Jordan algebras
  of dimension at most $4$ is distinguished,
  except in the case of the
  nilpotent Jordan algebra $\Ja_{4,63}$ which is not. 
\end{thm}

{\footnotesize
  \begin{table}
    \begin{center}
      \begin{tabular}{|c|c|c|}
      \hline
        \textsl{Isomorphism class} & \textsl{Multiplication table} &Properties\\
        \hline
$\Ja_{1,1}$ &$e_1^2=e_1$ & A, S \\
        \hline\end{tabular}
    \end{center}
    \caption{Complex $1$-dimensional Jordan algebras and their 
solitons.}\label{J1}
    \end{table}    
}

{\footnotesize
\begin{table}
    \begin{center}
      \begin{tabular}{|c|c|c|c|}
      \hline
      \textsl{$G$-orbits} & \textsl{Soliton type}  &$\beta$ & $E_1$ \\ 
      \hline 
      $\mathfrak A_{1,1}$ & $(0;1)$ & $\mathrm{diag}(-1)$ & $1$ \\
\hline
      \end{tabular}
    \end{center}
    \caption{The stratification of $\mathcal J_1$.}\label{J1bis}
    \end{table}
    }

{\footnotesize
  \begin{table}
    \begin{center}
      \begin{tabular}{|c|c|c|}
      \hline
        \textsl{Isomorphism class} & \textsl{Multiplication table} &
                                                                    Properties\\
        
      \hline 
        $\mathfrak A_{2,1}$ &  $e_1^2=e_1$, $e_1n_1=n_1$ & A, U \\
        $\mathfrak A_{2,2}$ &  $e_1^2=e_1$, $e_1n_1=\frac12n_1$ & - \\
        $\mathfrak A_{2,3}$ &  $n_1^2=n_2$ & A, N \\
        $\mathfrak A_{2,4}=(\Ja_{1,1})^2$ &  $e_1^2=e_1$, $e_2^2=e_2$ & A, SS, D \\
        $\mathfrak A_{2,5}=\Ja_{1,1}\times\Jt$ &  $e_1^2=e_1$, $n_1^2=0$ & A, D \\
        \hline\end{tabular}
    \end{center}
    \caption{Complex $2$-dimensional Jordan algebras and their solitons.}\label{J2}
    \end{table}    
}
{\footnotesize
   \begin{table}
    \begin{center}
      \begin{tabular}{|c|c|c|c|}
      \hline
      \textsl{$G$-orbits} & \textsl{Soliton type}  &$\beta$ & $E_2$ \\ 
      \hline 
      $\mathfrak A_{2,4}$ & $(0;2)$ & $\mathrm{diag}(-\tfrac12,-\tfrac12)$ & $1/2$ \\
      $\mathfrak A_{2,1}$, $\mathfrak A_{2,5}$, $\mathfrak A_{2,2}$ &$(0<1;1,1)$  &$\mathrm{diag}(-1,0)$ & $1$ \\
      $\mathfrak A_{2,3}$ & $(1<2;1,1)$ & $\mathrm{diag}(-2,1)$ & $5$\\
      \hline
      \end{tabular}
    \end{center}
    \caption{The stratification of $\mathcal J_2$.}\label{J2bis}
    \end{table}
   }
  {\footnotesize
  \begin{table}
\begin{center}
      \begin{tabular}{|c|c|c|}
      \hline
        \textsl{Isom class} & \textsl{Multiplication table} &
                                                                    Properties\\
        
      \hline 
        $\mathfrak A_{3,1}=(\Ja_{1,1})^3$ &  $e_1^2=e_1$, $e_2^2=e_2$, $e_3^2=e_3$ & SS, A, D\\
        $\mathfrak A_{3,2}$ &  $e_1^2=e_1$, $e_2^2=e_3^2=\frac{\sqrt5}2e_1$, $e_1e_2=e_2$, $e_1e_3=e_3$ & S\\
        $\mathfrak A_{3,3}=\Ja_{2,1}\times\Ja_{1,1}$ &  $e_1^2=e_1$, $e_2^2=\sqrt3e_2$, $e_1n_1=n_1$ & 
A, U, D\\
        $\mathfrak A_{3,4}$ &  $e_1^2=e_1$, $e_2^2=\sqrt{\frac53}e_1$, 
$e_1e_2=e_2$, $e_1n_1=n_1$ & U\\
        $\mathfrak A_{3,5}=\Ja_{2,2}\times\Ja_{1,1}$ &  $e_1^2=e_1$, $e_2^2=\sqrt{\frac32}e_2$, $e_1n_1=\frac12n_1$ & D\\
        $\mathfrak A_{3,6}=(\Ja_{1,1})^2\times\Jt$ &  $e_1^2=e_1$, $e_2^2=e_2$, $n_1^2=0$ & A, D\\
        $\mathfrak A_{3,7}$ & $e_1^2 = e_1$ , $e_1 n_1 = n_1$,
                         $e_1 n_2 = n_2$, $n_1^2 = n_2$ & A, U \\
        $\mathfrak A_{3,8}$ & $e_1^2 = e_1$ , $e_1 n_1 = n_1$,
                         $e_1 n_2 = n_2$& A, U\\
          $\mathfrak A_{3,9}=\Ja_{2,1}\times\Jt$ & $e_1^2 = e_1$ , $e_1 n_1 = n_1$, $n_2^2=0$
& A, D \\
        $\mathfrak A_{3,{10}}$ & $e_1^2 = e_1$ , $e_1 n_1 = \frac12n_1$,
                            $e_1 n_2 = n_2$, $n_1^2 = \sqrt{\frac7{10}}n_2$ & - \\
        $\mathfrak A_{3,{11}}$ & $e_1^2 = e_1$ , $e_1 n_1 = \frac12n_1$,
                            $e_1 n_2 = n_2$ & - \\
        $\mathfrak A_{3,{12}}$ & $e_1^2 = e_1$ , $e_1 n_1 = \frac12n_1$,
                            $e_1 n_2 = \frac12n_2$ & - \\
        $\mathfrak A_{3,{13}}$ & $e_1^2 = e_1$ , $e_1 n_1 = \frac12n_1$,
                            $n_1^2=\sqrt{\frac3{10}}n_2$ & - \\
        $\mathfrak A_{3,{14}}=\Ja_{2,2}\times\Jt$ & $e_1^2 = e_1$ , $e_1 n_1 = \frac12n_1$, $n_2^2=0$ & D \\
        $\mathfrak A_{3,{15}}=\Ja_{2,3}\times\Ja_{1,1}$ & $e_1^2=\sqrt 5e_1$, $n_1^2=n_2$ &  A, D \\
        $\mathfrak A_{3,{16}}=\Ja_{1,1}\times(\Jt)^2$ & $e_1^2=e_1$, $n_1^2=n_2^2=0$ &  A, D \\

        $\mathfrak A_{3,{17}}$ & $n_1^2=n_2$, $n_1n_2=n_3$ &  A, N \\
        $\mathfrak A_{3,{18}}$ & $n_1n_2=n_3$ &  A, N \\
         $\mathfrak A_{3,{19}}=\Ja_{2,3}\times\Jt$ & $n_1^2=n_2$, $n_3^2=0$ &  A, N, D \\
        \hline\end{tabular}
    \end{center}
    \caption{Complex $3$-dimensional Jordan algebras and their solitons.}\label{J3}
\end{table}
}

{\footnotesize
    \begin{table}
      \begin{center}
      \begin{tabular}{|c|c|c|c|}
      \hline
      \textsl{$G$-orbits} & \textsl{Soliton type}  &$\beta$ & $E_3$ \\ 
      \hline 
      $\mathfrak A_{3,1}$, $\mathfrak A_{3,2}$ & $(0;3)$ & $\mathrm{diag}(-\tfrac13,-\tfrac13,-\tfrac13)$ & $1/3$ \\
      $\mathfrak A_{3,3}$, $\mathfrak A_{3,4}$, $\mathfrak A_{3,5}$, $\mathfrak A_{3,6}$ & $(0<1;2,1)$  &$\mathrm{diag}(-\tfrac12,-\tfrac12,0)$ & $1/2$ \\
        $\mathfrak A_{3,7}$, $\mathfrak A_{3,{10}}$, $\mathfrak A_{3,{13}}$,
        $\mathfrak A_{3,{15}}$ & $(0<1<2;1,1,1)$ & $\mathrm{diag}(-\tfrac56,-\tfrac13,\tfrac16)$ & $5/6$\\
        $\mathfrak A_{3,8}$, $\mathfrak A_{3,9}$, $\mathfrak A_{3,{11}}$, $\mathfrak A_{3,{12}}$,
        $\mathfrak A_{3,{14}}$, $\mathfrak A_{3,{16}}$
                          & $(0<1;1,2)$ &  $\mathrm{diag}(-1,0,0)$ & $1$\\

        $\mathfrak A_{3,{17}}$
                          & $(1<2<3;1,1,1)$ &  $\mathrm{diag}(-\tfrac43,-\tfrac13,\tfrac23)$ & $7/3$\\
        $\mathfrak A_{3,{18}}$
                          & $(1<2;2,1)$ &  $\mathrm{diag}(-1,-1,1)$ & $3$\\
        $\mathfrak A_{3,{19}}$
                          & $(3<5<6;1,1,1)$ &  $\mathrm{diag}(-2,1,0)$ & $5$\\

      \hline
      \end{tabular}
    \end{center}
    \caption{The stratification of $\mathcal J_3$.}\label{J3bis}
    \end{table}
}

{\footnotesize
      \begin{table}
\begin{center}
      \begin{tabular}{|c|c|c|}
      \hline
        \textsl{Isom class} & \textsl{Multiplication table} &
                                                                    Properties\\
        
        \hline
$\Ja_{4,1}$ & $e_1^2=e_1$, $e_2^2=\frac{\sqrt5}2e_1$,
$e_3^2=\frac{\sqrt5}2e_1$, $e_1e_2=e_2$, $e_1e_3=e_3$, $e_4^2=\sqrt{\frac52}e_4$ 
& SS, D \\
$\Ja_{4,2}$ & $e_1^2=e_1$, $e_1e_2=e_2$, $e_1e_3=e_3$, $e_1e_4=e_4$,
$e_2e_3=\sqrt{\frac75}e_1$, $e_4^2=\sqrt{\frac75}e_1$ & S \\
$\Ja_{4,3}=(\Ja_{1,1})^4$ & $e_1^2=e_1$, $e_2^2=e_2$, $e_3^2=e_3$, $e_4^2=e_4$ &
SS, A, D \\
$\Ja_{4,4}=\Ja_{2,1}\times(\Ja_{1,1})^2$ & $e_1^2=e_1$, $e_1n_1=n_1$, $e_2^2=\sqrt3 e_2$, $e_3^2=\sqrt3 e_3$ & U, A, D\\
 $\Ja_{4,5}=(\Ja_{1,1})^3\times\mathfrak T$ & $e_1^2=e_1$, $e_2^2=e_2$, $e_3^2=e_3$, $n_1^2=0$ & A, D \\ 
$\Ja_{4,6}=\Ja_{2,2}\times(\Ja_{1,1})^2$ & $e_1^2=e_1$, $e_2^2=\sqrt{\frac32}e_2$,
$e_3^2=\sqrt{\frac32}e_3$, $e_1n_1=\frac12n_1$ & D \\
$\Ja_{4,7}=\Ja_{3,4}\times\Ja_{1,1}$ & 
$e_1^2=e_1$, $e_2^2=\sqrt{\frac53}e_1$, $e_1e_2=e_2$, $e_1n_1=n_1$,
$e_3^2=\sqrt{\frac{10}3}e_3$
& U, D \\ 
$\Ja_{4,8}=\Ja_{3,2}\times\mathfrak T$ & 
$e_1^2=e_1$, $e_2^2=\frac{\sqrt5}2e_1$, $e_3^2=\frac{\sqrt5}2e_1$,
$e_1e_2=e_2$, $e_1e_3=e_3$, $n_1^2=0$ 
& D \\
$\Ja_{4,9}$ & $e_1^2=e_1$, $e_2^2=\frac{\sqrt7}2e_1$, 
$e_3^2=\frac{\sqrt7}2e_1$, $e_1e_2=e_2$, $e_1e_3=e_3$, $e_1n_1=n_1$
& U\\
$\Ja_{4,10}=\Ja_{2,2}\times\Ja_{1,1}\times\Jt$ & $e_1^2=e_1$, $e_1n_1=\frac12n_1$,
$e_2^2=\sqrt{\frac32}e_2$, $n_2^2=0$ & D \\
$\Ja_{4,11}=\Ja_{3,4}\times\mathfrak T$ & $e_1^2=e_1$,
$e_2^2=\sqrt{\frac53}e_1$, $e_1e_2=e_2$, $e_1n_1=n_1$, $n_2^2=0$ & D \\
$\Ja_{4,12}=\Ja_{3,12}\times\Ja_{1,1}$ & $e_1^2=e_1$, $e_2^2=\sqrt2e_2$, $e_1n_1=\frac12n_1$, $e_1n_2=\frac12n_2$ & D\\
$\Ja_{4,13}=(\Ja_{2,2})^2$ & $e_1^2=e_1$, $e_2^2=e_2$, $e_1n_1=\frac12n_1$, $e_2n_2=\frac12n_2$ & D\\
$\Ja_{4,14}=\Ja_{3,11}\times\Ja_{1,1}$ & $e_1^2=e_1$, $e_1n_1=\frac12n_1$, $e_1n_2=n_2$,
$e_2^2=\sqrt{\frac72}e_2$ & D \\
        $\Ja_{4,15}=\Ja_{2,1}\times\Ja_{2,2}$ & $e_1^2=e_1$, $e_2^2=\sqrt2e_2$, $e_1n_1=n_1$, $e_2n_2=\frac1{\sqrt2}n_2$ & D \\
\multirow{5}{*}{$\Ja_{4,16}$} & \parbox[t]{3in}{$e_1^2=k(\cos^3 t-\sin^3t)e_1+k^2\cos t\sin t (\cos t + \sin t)e_2$, $e_2^2=k^{-1}\cos t\sin t(-\cos t+\sin t)e_1+(\cos^3 t + \sin^3 t)e_2$,
$e_1e_2=\cos t\sin t((\cos t+\sin t)e_1 + k (-\cos t+\sin t)e_2)$,
$e_1n_1=\frac k2(\cos t-\sin t)n_1$, $e_1n_2=\frac k2\cos t n_2$,\\
$e_2n_1=\frac12(\cos t + \sin t)n_1$, $e_2n_2=\frac12\sin t n_2$} &
\multirow{5}{*}{\parbox[t]{.8in}{$k\approx1.20577$, $t\approx1.22166$}} \\
\multirow{5}{*}{$\Ja_{4,17}$} & \parbox[t]{3in}{$e_1^2=k(\cos^3 t-\sin^3t)e_1+k^2\cos t\sin t (\cos t + \sin t)e_2$, $e_2^2=k^{-1}\cos t\sin t(-\cos t+\sin t)e_1+(\cos^3 t + \sin^3 t)e_2$,
$e_1e_2=\cos t\sin t(\cos t+\sin t)e_1 + k \cos t\sin t(-\cos t+\sin t)e_2$,
$e_1n_1=\frac k2(\cos t-\sin t)n_1$, $e_1n_2=k\cos t n_2$,
$e_2n_1=\frac12(\cos t + \sin t)n_1$, $e_2n_2=\sin t n_2$} &
\multirow{5}{*}{\parbox[t]{.8in}{$k\approx1.54492$, $t\approx1.45358$, U}} \\
        $\Ja_{4,18}$ & $e_1^2=e_1$, $e_2^2=\sqrt{\frac73}e_1$, $e_1e_2=e_2$,
 $e_1n_1=n_1$, $e_1n_2=n_2$     & U \\
$\Ja_{4,19}=(\Ja_{1,1})^2\times(\mathfrak T)^2$ & $e_1^2=e_1$, $e_2^2=e_2$, $n_1^2=n_2^2=0$ & A, D\\
$\Ja_{4,20}=\Ja_{2,1}\times\Ja_{1,1}\times\Jt$ & $e_1^2=e_1$, $e_2^2=\sqrt3 e_2$, $e_1n_1=n_1$, $n_2^2=0$ & A, D \\
$\Ja_{4,21}=\Ja_{3,8}\times\Ja_{1,1}$ & $e_1^2=e_1$, $e_2^2=\sqrt5e_2$, $e_1n_1=n_1$, $e_1n_2=n_2$ & A, U, D\\
$\Ja_{4,22}=(\Ja_{2,1})^2$ & $e_1^2=e_1$, $e_2^2=e_2$, $e_1n_1=n_1$, $e_2n_2=n_2$ &
                                                                                    A, U, D \\
$\Ja_{4,23}=\Ja_{3,13}\times\Ja_{1,1}$ & $e_1^2=e_1$, $e_1n_1=\frac12n_1$,
$n_1^2=\sqrt{\frac3{10}}n_2$, $e_2^2=\sqrt{\frac32}e_2$ & D \\
$\Ja_{4,24}=\Ja_{3,10}\times\Ja_{1,1}$ & $e_1^2=e_1$, $e_1n_1=\frac12n_1$,
$e_1n_2=n_2$, $n_1^2=\sqrt{\frac{7}{10}}n_2$, $e_2^2=\sqrt{\frac72}e_2$ & 
D \\
\multirow{5}{*}{$\Ja_{4,25}$} & \parbox[t]{3in}{$e_1^2=k(\cos^3 t-\sin^3t)e_1+k^2\cos t\sin t (\cos t + \sin t)e_2$, $e_2^2=k^{-1}\cos t\sin t(-\cos t+\sin t)e_1+(\cos^3 t + \sin^3 t)e_2$,
$e_1e_2=\cos t\sin t(\cos t+\sin t)e_1 + k \cos t\sin t(-\cos t+\sin t)e_2$,
$e_1n_1=\frac k2(\cos t-\sin t)n_1$, $e_1n_2=k\cos t n_2$,
$e_2n_1=\frac12(\cos t + \sin t)n_1$, $e_2n_2=\sin t n_2$,
$n_1^2=\ell n_2$} &
\multirow{5}{*}{\parbox[t]{.8in}{$k\approx1.54492$, $t\approx1.45358$, $\ell\approx0.836502$, U}} \\
$\Ja_{4,26}=\Ja_{2,3}\times(\Ja_{1,1})^2$ & $e_1^2=e_1$, $e_2^2=e_2$, 
$n_1^2=\frac1{\sqrt5}n_2$, $n_2^2=0$ & A, D \\
$\Ja_{4,27}=\Ja_{3,7}\times\Ja_{1,1}$ & $e_1^2=e_1$, $e_1n_1=n_1$, $e_1n_2=n_2$,
$n_1^2=n_2$, $e_2^2=\sqrt5e_2$ & U, A, D \\
$\Ja_{4,28}=\Ja_{2,2}\times(\mathfrak T)^2$ & $e_1^2=e_1$, $e_1n_1=\frac12n_1$,
$n_2^2=n_3^2=0$ & D \\
$\Ja_{4,29}=\Ja_{3,11}\times\mathfrak T$ & $e_1^2=e_1$, $e_1n_1=\frac12n_1$,
$e_1n_2=n_2$, $n_3^2=0$ & D \\
$\Ja_{4,30}=\Ja_{3,12}\times\mathfrak T$ & $e_1^2=e_1$, $e_1n_1=\frac12n_1$,
$e_1n_2=\frac12n_2$, $n_3^2=0$ & D \\  
$\Ja_{4,31}$ & $e_1^2=e_1$, $e_1n_1=n_1$, $e_1n_2=n_2$, $e_1n_3=\frac12n_3$ 
& - \\
$\Ja_{4,32}$ & $e_1^2=e_1$, $e_1n_1=\frac12n_1$, $e_1n_2=\frac12n_2$, $e_1n_3=n_3$& - \\
$\Ja_{4,33}$ & $e_1^2=e_1$, $e_1n_1=\frac12n_1$, $e_1n_2=\frac12n_2$, $e_1n_3=\frac12n_3$& - \\
$\Ja_{4,34}=\Ja_{1,1}\times(\Jt)^3$ & $e_1^2=e_1$, $n_1^2=n_2^2=n_3^2=0$ & A, D \\
$\Ja_{4,35}=\Ja_{2,1}\times(\Jt)^2$ & $e_1^2=e_1$, $e_1n_1=n_1$, $n_2^2=n_3^2=0$
& A, D \\
$\Ja_{4,36}$ & $e_1^2=e_1$, $e_1n_1=n_1$, $e_1n_2=n_2$, $e_1n_3=n_3$ & U, A \\
$\Ja_{4,37}=\Ja_{3,8}\times\Jt$ & $e_1^2=e_1$, $e_1n_1=n_1$, $e_1n_2=n_2$, 
$n_3^2=0$ & A, D \\
$\Ja_{4,38}=\Ja_{3,17}\times\Ja_{1,1}$ & $e_1^2=\sqrt7e_1$, $n_1^2=n_2$, $n_1n_2=n_3$ & A, D \\
$\Ja_{4,39}$ &   $e_1^2=e_1$, $e_1n_1=n_1$, $e_1n_2=n_2$, $e_1n_3=n_3$,
$n_1^2=n_2$, $n_1n_2=n_3$ & U, A \\
$\Ja_{4,40}=\Ja_{2,3}\times\Ja_{1,1}\times\Jt$ & $n_1^2=n_2$ , $e_1^2=\sqrt5e_1$, $n_3^2=0$ 
& D \\
      \hline\end{tabular}
    \end{center}
    \caption{Complex $4$-dimensional Jordan algebras and 
their solitons, first part.}\label{J4i}
\end{table}

      \begin{table}
\begin{center}
      \begin{tabular}{|c|c|c|}
      \hline
        \textsl{Isom class} & \textsl{Multiplication table} &
                                                                    Properties\\
        
        \hline
$\Ja_{4,41}=\Ja_{3,18}\times\Ja_{1,1}$ & $e_1^2=\sqrt6e_1$, $n_1n_2=n_3$
& A, D \\
$\Ja_{4,42}$ & $e_1^2=e_1$, $e_1n_1=n_1$, $e_1n_2=n_2$, $e_1n_3=n_3$,
$n_1^2=\sqrt{\frac75}n_2$ & U, A \\ 
$\Ja_{4,43}$ & $e_1^2=e_1$, $e_1n_1=n_1$, $e_1n_2=n_2$, $e_1n_3=n_3$,
$n_1^2=n_2^2=\sqrt{\frac76}n_3$ & U, A \\
$\Ja_{4,44}=\Ja_{3,13}\times\Jt$ & $e_1^2=\sqrt{\frac{10}3}e_1$, 
$e_1n_1=\sqrt{\frac56}n_1$, $n_1^2=n_2$, $n_3^2=0$ & D \\
$\Ja_{4,45}$ & $e_1^2=2e_1$, $e_1n_1=n_1$, $n_1^2=n_2^2=n_3$ & - \\
$\Ja_{4,46}=\Ja_{2,2}\times\Ja_{2,3}$ & $e_1^2=e_1$, $e_1n_1=\frac12n_1$,
$n_2^2=\sqrt{\frac3{10}}n_3$ & D \\
$\Ja_{4,47}=\Ja_{2,1}\times\Ja_{2,3}$ & $e_1^2=e_1$, $e_1n_1=n_1$,
$n_2^2=\sqrt{\frac35}n_3$ & A, D \\
$\Ja_{4,48}$ & $e_1^2=e_1$, $e_1n_1=\frac12n_1$, $e_1n_2=\frac12n_2$,
$n_1^2=\sqrt{\frac25}n_3$ & - \\
$\Ja_{4,49}$ & $e_1^2=e_1$, $e_1n_1=\frac12n_1$, $e_1n_2=\frac12n_2$,
$n_1^2=n_2^2=\frac1{\sqrt3}n_3$ & - \\
$\Ja_{4,50}$ & $e_1^2=e_1$, $e_1n_2=\frac12n_2$, $e_1n_3=\frac12n_3$,
$n_1n_2=\frac1{\sqrt3}n_3$ & - \\
$\Ja_{4,51}=\Ja_{3,10}\times\Jt$ & $e_1^2=e_1$, $e_1n_1=\frac12n_1$, $e_1n_2=n_2$,
$n_1^2=\sqrt{\frac7{10}}n_2$, $n_3^2=0$ & D \\
$\Ja_{4,52}$ & $e_1^2=e_1$, $e_1n_1=\frac12n_1$, $e_1n_2=n_2$, 
$n_1^2=\sqrt{\frac7{10}}n_3$ & - \\
$\Ja_{4,53}$ & $e_1^2=e_1$, $e_1n_1=\frac12n_1$, 
$e_1n_2=\frac12(n_2+\alpha n_3)=\alpha e_1n_3$,
$n_1^2=\beta n_3$ & $\alpha^2=\frac{\sqrt{345}-5}{20}$, $\beta^2=\frac{\sqrt{345}+45}{80}$ \\
$\Ja_{4,54}=\Ja_{3,7}\times\Jt$ & $e_1^2=e_1$, $e_1n_1=n_1$, 
$e_1n_2=n_2$, $n_1^2=n_2$, $n_3^2=0$ & A, D \\
$\Ja_{4,55}$ &  $e_1^2=e_1$, $e_1n_1=n_1$, $e_1n_2=n_2$, $e_1n_3=\frac12 n_3$,
$n_3^2=\sqrt{\frac{11}{10}}n_2$ & - \\
$\Ja_{4,56}$ &  $e_1^2=e_1$, $e_1n_1=n_1$, $e_1n_2=n_2$, $e_1n_3=\frac12 n_3$,
$n_1^2=\sqrt{\frac{11}{10}}n_2$ & - \\
$\Ja_{4,57}$ &  $e_1^2=e_1$, $e_1n_1=n_1$, $e_1n_2=n_2$, $e_1n_3=\frac12 n_3$,
$n_1^2=n_3^2=\sqrt{\frac{11}{12}}n_2$ & - \\
$\Ja_{4,58}$ &  $e_1^2=e_1$, $e_1n_1=n_1$, $e_1n_2=\frac12n_2$, $e_1n_3=\frac12 n_3$, $n_3^2=\frac{2}{\sqrt5}n_1$ & - \\
$\Ja_{4,59}$ &  $e_1^2=e_1$, $e_1n_1=n_1$, $e_1n_2=\frac12n_2$, $e_1n_3=\frac12 n_3$, $n_2^2=n_3^2=\sqrt{\frac23}n_1$ & - \\
$\Ja_{4,60}$ &  $e_1^2=e_1$, $e_1n_1=n_1$, $e_1n_2=\frac12n_2$, $e_1n_3=\frac12 n_3$, $n_1n_2=\sqrt{\frac23}n_3$ & - \\
$\Ja_{4,61}$ & $n_1^2=n_2$, $n_2^2=n_4$, $n_1n_2=n_3$, $n_1n_3=n_4$ & A, N \\
$\Ja_{4,62}$ & $n_1^2=n_2$, $n_4^2=2n_2$, $n_1n_2=\sqrt3n_3$ & N \\
$\Ja_{4,63}$ & $n_1n_2=n_3$, $n_1n_3=n_4$, $n_2^2=n_4$ & N, \mbox{not distinguished} \\
$\Ja_{4,64}$ & $n_1n_2=n_3$, $n_1n_3=n_4$ & N \\
$\Ja_{4,65}$ & $n_1^2= \frac2{\sqrt3}n_2$, $n_2n_3=n_4$ & N \\
$\Ja_{4,66}$ & $n_1^2=n_2$, $n_3^2=n_4$, $n_1n_2=\frac{\sqrt3}2n_4$ & A, N \\
$\Ja_{4,67}=\Ja_{3,17}\times\Jt$ & $n_1^2=n_2$, $n_1n_2=n_3$, $n_4^2=0$ 
& A, N, D \\
$\Ja_{4,68}=(\Ja_{2,3})^2$ & $n_1^2=n_2$, $n_3^2=n_4$ & A, N, D \\
$\Ja_{4,69}$ & $n_1^2=n_2$, $n_1n_3=\sqrt{\frac32}n_4$ & A, N\\
$\Ja_{4,70}$ & $n_1^2=n_3n_4=n_2$ & A, N\\
$\Ja_{4,71}=\Ja_{3,18}\times\Jt$ & $n_1n_2=n_3$, $n_4^2=0$ & A, N \\
$\Ja_{4,72}=\Ja_{2,3}\times(\Jt)^2$ & $n_1^2=n_2$, $n_3^2=n_4^2=0$ & A, N \\
  \hline\end{tabular}
    \end{center}
    \caption{Complex $4$-dimensional Jordan algebras and their
solitons, second part.}\label{J4ii}
\end{table}
}

{\small
    \begin{table}
      \begin{center}
      \begin{tabular}{|c|c|c|c|}
      \hline
      \textsl{$G$-orbits $\Ja_{4,k}$} & \textsl{Soliton type}  &$\beta$ & $E_4$ \\ 
\hline
  $1\leq k\leq3$ & $(0;4)$ & $\mathrm{diag}(-\tfrac14,-\tfrac14,-\tfrac14,-\tfrac14)$ & $1/4$ \\
  $4\leq k\leq9$  & $(0<1;3,1)$ & $\mathrm{diag}(-\tfrac13,-\tfrac13,-\tfrac13,0)$ & $1/3$ \\
 $23\leq k\leq27$  & $(0<1<2;2,1,1)$ & $\mathrm{diag}(-\tfrac5{11},-\tfrac5{11},-\tfrac2{11},\tfrac1{11})$ & $5/11$ \\
 $10\leq k\leq22$  & $(0<1;2,2)$ & $\mathrm{diag}(-\tfrac12,-\tfrac12,0,0)$ & $1/2$ \\
$38\leq k\leq39$  & $(0<1<2<3;1,1,1,1)$ & $\mathrm{diag}(-\tfrac7{10},-\tfrac4{10},-\tfrac1{10},\tfrac2{10})$ & $7/10$ \\
$k\in\{41,43,45,49,50,57,59,60\}$ & $(0<1<2;1,2,1)$ & $\mathrm{diag}(-\tfrac34,-\tfrac14,-\tfrac14,\tfrac14)$ & $3/4$ \\
$k=53$ & $(0<1<2;1,1,2)$ & $\mathrm{diag}(-\tfrac9{11},-\tfrac4{11},\tfrac1{11},\tfrac1{11})$ & 
$9/11$ \\ 
$k\in\{40,42,44,46,47,48,51,52,54,55,56,58\}$  
& $(0<3<5<6;1,1,1,1)$ & $\mathrm{diag}(-\tfrac56,-\tfrac13,0,\tfrac16)$ & 
$5/6$ \\
$28\leq k\leq37$  & $(0<1;1,3)$ & $\mathrm{diag}(-1,0,0,0)$ & $1$ \\
$k=62$  & $(1<2<3;2,1,1)$ & $\mathrm{diag}(-\tfrac8{11},-\tfrac8{11},-\tfrac1{11},\tfrac6{11})$ & $15/11$ \\
$k=65$  & $(3<4<6<10;1,1,1,1)$ & $\mathrm{diag}(-\tfrac45,-\tfrac35,-\tfrac15,\tfrac35)$ & $7/5$ \\
$k=61$, $63$, $64$  & $(1<2<3<4;1,1,1,1)$ & $\mathrm{diag}(-1,-1/2,0,1/2)$ & $3/2$ \\
        $k=66$  & $(2<3<4<6;1,1,1,1)$ & $\mathrm{diag}(-1,-\tfrac17,-\tfrac47,\tfrac57)$ & $13/7$ \\
        $k=70$  & $(1<2;3,1)$ & $\mathrm{diag}(-\tfrac23,-\tfrac23,-\tfrac23,1)$ & $7/3$ \\
         $k=67$  & $(3<6<7<9;1,1,1,1)$ & $\mathrm{diag}(-\tfrac43,-\tfrac13,0,\tfrac23)$ & $7/3$ \\
$k=68$  & $(1<2;2,2)$ & $\mathrm{diag}(-1,-1,\tfrac12,\tfrac12)$ & $5/2$ \\
$k=69$  & $(3<4<6<7;1,1,1,1)$ & $\mathrm{diag}(-\tfrac54,-\tfrac34,\tfrac14,\tfrac34)$ & $11/4$ \\
$k=71$  & $(2<3<4;2,1,1)$ & $\mathrm{diag}(-1,-1,0,1)$ & $3$ \\
$k=72$  & $(3<5<6;1,2,1)$ & $\mathrm{diag}(-2,0,0,1)$ & $5$ \\
\hline
 \end{tabular}
    \end{center}
    \caption{The stratification of $\mathcal J_4$.}\label{J4bis}
    \end{table}
}

In Tables~\ref{J1}, \ref{J2}, \ref{J3}, \ref{J4i} and 
\ref{J4ii} we list the complex Jordan solitons of dimension at most $4$.
The coefficients in the multiplication table are chosen so that 
the algebra in each isomorphism type is a soliton, except 
$\Ja_{4,63}$ for which there is no soliton. In the column of properties
we indicate whether the algebra is associative (A), simple (S), 
semisimple (SS), nilpotent (N), unitary (U), or decomposable (D); the absence 
of such a qualification means its negation, except that a simple Jordan
algebra is semisimple, and a semisimple Jordan algebra is unitary. 
In the decomposable case we exhibit the minimal decomposition; the notation
$\Jt$ refers to the trivial (one-dimensional) algebra.  
In Tables~\ref{J1bis}, \ref{J2bis}, \ref{J3bis}, \ref{J4bis} 
we describe the strata of $\mathcal J_n$ for $n\leq4$: 
the orbits they contain,  the soliton type, 
the parameter ($\beta$), and the value of the energy ($E_n=||\beta||^2$). 

The property of a $G$- orbit to be distinguished is independent of the
choice of Hermitean product~$Q$. Indeed there is an action of $G$
on the space of Hermitean products on $\C^n$ given by
\[ g\cdot Q(a,b) = Q(g^{-1}a,g^{-1}b)) \]
for $g\in G$ and $a$, $b\in\C^n$. This action is transitive, so,
if we fix $Q$, any other
Hermitean product on $\C^n$ is of the form $g\cdot Q$ for some $g\in G$.
Using a superscript to denote the moment matrices associated to different
Hermitean products, it is easy to see that
\begin{equation}\label{change-Q}
  M^{g\cdot Q}_\mu= gM^Q_{g^{-1}\cdot\mu}g^{-1}. 
\end{equation}
Hence the claim follows from the structure result
for solitons (Proposition~\ref{prop:struct}(a)).

Thinking in a dual way, fix $Q$ and suppose the isomorphism class of a Jordan
algebra $\Ja=(\C^n,\mu)$ is distinguished. Then $g\cdot\mu$ is a soliton for some
$g\in G$, and the relation~(\ref{change-Q}) can be used to see that 
$\mu$ is a soliton with respect to $g^{-1}\cdot Q$, that is, 
$g^{-1}\cdot Q$ is a 'best' metric on $\Ja$. It follows from Kirwan-Ness 
theory that the critical set of $E_n$ on a $G$-orbit, if non-empty, 
is a single $K$-orbit, so the best metric on $\Ja$, if existing, is 
unique. 

We point out that several formulae and results herein
are similar to those for Lie 
algebras~\cite{lauret2003}. An important difference of Jordan algebras from
Lie algebras is that, in general, left multiplications are not derivations. 
So arguments for Lie algebras that depend on this property will not work
in the setting of Jordan algebras. Of course this should not be seen as 
a problem, but as a central feature of Jordan algebras; and the real problem is 
knowing how to use this feature to our advantage.

The authors wish to thank I. Shestakov for stimulating discussions. 

\section{Preliminaries}

\subsection{Semisimplicity and the radical of Jordan algebras}\label{sec:ss}

 Most of the results in this subsection are due to Albert.
A good reference is~\cite[chap.~IV]{schafer}. 

A complex Jordan algebra $\mathfrak A$ is called
  \emph{simple} if $\mathfrak A^2\neq0$ and $\mathfrak A$
  have no nontrivial ideals, and it is called
  \emph{semisimple} if it is a direct product of
  simple Jordan algebras. A semisimple Jordan algebra 
  has an identity element. 

Let $\mathfrak A$ be a complex Jordan algebra. 
An element $x\in\Ja$ is called \emph{nilpotent} if $x^k=0$ for some
integer $k\geq2$; this is equivalent to the left-multiplication 
operator $L_x$ being a nilpotent operator. There
is an ideal $\mathrm{Rad}(\Ja)$, called the \emph{radical} of 
$\Ja$, which is the unique maximal
nilideal of $\Ja$ (that is, the maximal ideal consisting entirely of 
nilpotent elements). Furthermore, $\mathrm{Rad(\Ja)}$ 
is nilpotent in the sense that there
is an integer $t$ with the property that any product 
$x_1\cdots x_t$  of $t$ elements
from $\mathrm{Rad}(\Ja)$ is zero; 
hence $\mathrm{Rad}(\Ja)$ 
is also the unique maximal nilpotent ideal of $\Ja$. 
With this definition of radical, the quotient algebra 
$\Ja/\mathrm{Rad}(\Ja)$ is always semisimple, and 
$\Ja$ is semisimple if and only if $\mathrm{Rad}(\Ja)=0$.  

A result of Albert characterizes the radical $\mathrm{Rad}(\Ja)$ as the 
kernel of the symmetric bilinear form 
$\tau:\mathfrak A\times\mathfrak A\to\C$ defined 
  by $\tau(x,y)=\mathrm{Tr}(L_{xy})$ for $x$, $y\in\mathfrak A$.
 
The Wedderburn Principal Theorem for Jordan algebras, 
proved by Albert and Penico~\cite{penico}, states that  
any complex Jordan algebra $\mathfrak A$ can be written
as a vector space direct sum  $\mathfrak A=\mathfrak S+\mathfrak N$ 
for some maximal semisimple subalgebra $\mathfrak S$ of~$\mathfrak A$
isomorphic to $\Ja/\Jn$, where 
$\mathfrak N=\mathrm{Rad}(\mathfrak A)$.

\subsection{A review of GIT}

Let $\pi:G\to GL(V)$ be a rational representation of a connected
complex linearly reductive group $G$ on a finite-dimensional complex
vector space $V$. Fix a maximal compact subgroup $K$ of $G$ (necessarily
connected and unique up to inner automorphism of $G$),
and a $K$-invariant Hermitian inner product $\langle,\rangle$
on $V$. Then the Lie algebra $\Lk$ of $K$ is a real form of the Lie algebra
$\Lg$ of $G$, namely, $\Lg=\Lk+i\Lk$. We also fix an $\Ad_K$-invariant
Hermitian inner product on $\Lg$, denoted $(,)$, 
which is positive definite (resp. negative definite) on $i\Lk$ (resp.~$\Lk$). 

The \emph{moment map} is the map $m:V\setminus\{0\}\to\Lg$ defined by
\begin{equation}\label{moment-map}
 (m(v),X)= \frac1{2||v||^2}\frac{d}{dt}\Big|_{t=0}||e^{tX}\cdot v||^2=\frac{\langle X\cdot v,v\rangle}{\langle v,v\rangle}, 
\end{equation}
where we write $g\cdot v:=\pi(g)v$ and $X\cdot v:=d\pi(X)v$, for all $g\in G$, 
$X\in\Lg$, $v\in V\setminus\{0\}$. We have that 
$m(v)$ controls the norm of vectors in a neighborhood of $v$
in the orbit $G\cdot v$ orthogonal to~$v$, and $-m(v)v$ is the direction of fastest
decrease in the norm. 
Note that $m(v)\in i\Lk$ and $m$ is $K$-invariant by 
$K$-invariance of $\langle,\rangle$. 
A \emph{minimal} vector is vector of minimal length in its orbit.  
The Kempf-Ness theorem says that minimal vectors correspond to the
zeros of the moment map, and  an orbit is closed if and only
it contains a minimal vector;
further the closure of any $G$-orbit in $V$ contains exactly one $K$-orbit of
minimal vectors~\cite{kempf-ness}.
The square norm of $m$ yields the energy functional
\[ E: V\setminus\{0\}\to\R,\qquad E(v)=||m(v)||^2. \]
Of course $E$ is $K$-invariant and scale-invariant.
It follows from the above that  a $G$-orbit in $V$ is 
closed if and only if it meets~$E^{-1}(0)$.

On the other hand, the union of all the orbits containing the origin in its
closure comprise the so-called \emph{null cone} $N$ of $V$. 
The Kirwan-Ness theorem says that the non-minimal critical points 
of $E$, i.~e.~those critical points $v$ of $E$ with $E(v)>0$,
all occur in the null cone $N$ and determine a stratification of 
$N\setminus\{0\}$~\cite{kirwan,ness}. 
So there is a set of distinguished orbits in $N$, namely those 
containing a critical point of $E$, which play a role similar to the 
closed orbits in $V\setminus N$.

Since $(dm_v(\xi),X)=2\Re\langle X\cdot  v,\xi\rangle/||v||^2$
for $X\in i\Lk$, $\xi\in T_vV=V$, $\xi\perp v$,
we easily obtain for the differential of
$E$ at a point $v\in V\setminus\{0\}$,
\[ dE_v(\xi)=2\Re(dm_v(\xi),m(v)) = \frac{4\Re\langle m(v)\cdot v,\xi\rangle}{||v||^2}. \]
Now $E$ is constant in the complex radial direction, and
the component of $m(v)\cdot v$ in the direction of $v$ is
$\frac{\langle m(v)\cdot v,v\rangle}{||v||^2}v=||m(v)||^2v$,
so the gradient of $E$ at~$v$ is given by 
\begin{equation}\label{grad-e}
 \nabla E(v) = \frac4{||v||^2}(m(v)\cdot v-||m(v)||^2v). 
\end{equation}
Therefore $v\in V\setminus\{0\}$ is a critical point of $E$ if and only if
\begin{equation}\label{self-similar}
 m(v)\cdot v\in\R v. 
\end{equation}
Because of this self-similarity characteristic, 
a critical point  of $E$ is called a \emph{soliton}
(following~\cite{lauret2020}).
We will also say a $G$-orbit  containing 
a soliton is \emph{distinguished}. Of course every minimal vector is a 
soliton. For the convenience of the reader,
we next collect the results of Kirwan-Ness theory related to solitons 
that we shall need. 

\begin{thm}[Kirwan-Ness]\label{kirwan-ness}
With the above notation, we have:
\begin{enumerate}
\item The subset of solitons in a given $G$-orbit is either empty
  or consists of precisely one $K$-orbit, up to scaling. 
\item Every soliton $v$ is a minimum of $E$ on $G\cdot v$. Thus solitons
  are the vectors closest to being minimal in their $G$-orbits.
\item The solitons which are not minimal vectors all occur in the 
null cone $N$ of $V$.
\item The flow of $-\nabla E$ starting at $v$ stays in $G\cdot v$
  and converges to a soliton $w\in\overline{G\cdot v}$ as $t\to\pm\infty$.
  There is precisely one $K$-orbit up to scaling
  of solitons $z\in\overline{G\cdot v}$
  such that $m(z)\in\Ad_Km(w)$, which consists of the limit-set
  of $G\cdot v$.
\item The critical set of $E$ is a finite disjoint union of closed
  subsets $\{C_\beta\}_{\beta\in\mathcal B}$
  indexed by a finite set $\mathcal B$ of adjoint $K$-orbits in $\Lk$
  (or points in a positive Weyl chamber). The corresponding
  stable manifolds $\{S_\beta\}_{\beta\in\mathcal B}$ form a
  finite $G$-invariant stratification of $N\setminus\{0\}$ by 
locally closed irreducible nonsingular subvarieties. 
\end{enumerate}
\end{thm}

\section{The moment map for commutative algebras}

In this section we specialize to the case of $n$-dimensional commutative 
algebras and, in particular, Jordan algebras. 
Let $V_n=S^2(\C^{n*})\otimes\C^n$ be the space of symmetric bilinear maps 
$\C^n\times\C^n\to\C^n$, which we identify with 
(non-necessarily associative) commutative algebras 
of dimension~$n$. There is a natural action of $G=GL(n,\C)$ 
on $V_n$ given by~(\ref{g-action}). 
The action of the Lie algebra
$\Lg=\mathfrak{gl}(n,\C)$ of $G$ on $V_n$,
obtained from differentiation of~(\ref{g-action}), is given by
\begin{equation}\label{A-action}
 A\cdot\mu(x,y) =A(\mu(x,y))-\mu(Ax,y)-\mu(x,Ay) 
\end{equation}
for all $A\in\Lg$, $\mu\in V_n$, and $x$, $y\in\C^n$. 
Note that the isotropy algebra of $\mu$, namely, the subalgebra of 
$\Lg$ consisting of elements $A\in\Lg$ satisfying $A\cdot\mu=0$,
is isomorphic to the Lie algebra $\mathrm{Der}(\mu)$ 
of derivations of $\mu$. 

Consider the canonical Hermitean product $\langle,\rangle$ of $\C^n$. 
The unitary group $K=U(n)$ is a maximal compact subgroup of $G$
and its Lie algebra $\Lk$ is a real form of $\Lg$. The Hermitean 
product canonically extends to a Hermitean product on $V_n$, denoted 
by the same symbol, namely, 
\[ \langle\mu,\nu\rangle =\sum_{ijk}\langle\mu(x_i,x_j),x_k\rangle
\overline{\langle\nu(x_i,x_j),x_k\rangle}, \]
where $x_1,\ldots,x_n$ is any fixed orthonormal basis of $\C^n$.
The elements of $K$ act on $V_n$ by unitary transformations, and 
those of $\Lk$ (resp.~$i\Lk$) 
act on $V_n$ by skew-Hermitean (resp. Hermitean) endomorphisms.
We also consider the $\Ad_K$-invariant Hermitean product on 
$\Lg$ given by $(A,B)=\mathrm{Tr}(AB^*)$, where~$A$, $B\in\Lg$.  

The moment map $m:V_n\setminus\{0\}\to i\Lk$ 
is defined as in~(\ref{moment-map}) and has the form
$m(\mu)=\frac1{||\mu||^2}M_\mu$, where $M_\mu$ is the moment matrix
of~$\mu$. 
We will obtain an explicit formula for $M_\mu$
in terms of the algebra structure of $\mu$ and the Hermitean
product momentarily. 
It turns out to be the same formula as that 
in~\cite{lauret2003} for skew-symmetric algebras.

Recall that $(M_\mu,A)=\langle A \cdot\mu,\mu\rangle$
for all~$A\in i\Lk$.

\begin{lem}\label{lem:moment}
For all $\mu\in V_n$ and $D\in\mathrm{Der}(\mu)$, we have:
\begin{enumerate}
\item[(a)] $\mathrm{Tr}(M_\mu)=-||\mu||^2$;
\item[(b)] $(M_\mu,D)=0$.
\end{enumerate}
\end{lem}

\Pf (a) Note that $\mathrm{Tr}(M_\mu)=(M_\mu,I)=\langle I\cdot\mu,\mu\rangle=-||\mu||^2$, as $I\cdot\mu=-\mu$. 

(b) We have  $(M_\mu,D)=\langle D\cdot\mu,\mu\rangle =0$, since
$D\cdot\mu=0$ by~(\ref{A-action}) and the derivation property. \EPf

\medskip

For $\mu\in V_n$, let $L^\mu_x:\C^n\to\C^n$
denote the left multiplication by $x\in\C^n$ in the algebra $\mu$,
that is, $L^\mu_xy=\mu(x,y)$ for all $y\in\C^n$. We also use a
superscript $()^*$ to denote the adjoint of a map $\C^n\to\C^n$ 
with respect the fixed Hermitean product. 

\begin{prop}\label{moment-jordan}
  The moment matrix of $\mu\in V_n\setminus\{0\}$ is given by  
\[ M_\mu = -2\sum_{i=1}^n L^{\mu*}_{x_i}L^\mu_{x_i}+\sum_{i=1}^nL^{\mu}_{x_i}L^{\mu*}_{x_i}, \]
where $(x_1,\ldots,x_n)$ is a orthonormal basis of $\C^n$. 
\end{prop}

\Pf Let $\mu\in V_n$ and $A\in i\Lk$. We have 
\begin{align*}
(M_\mu,A)&=\langle A\cdot\mu,\mu\rangle \\
         &=\sum_{ijs}\langle (A\mu)(x_s,x_i),x_j\rangle\overline{\langle \mu(x_s,x_i),x_j\rangle}\\
&=\sum_{ijs}
\langle A(\mu(x_s,x_i))-\mu(Ax_s,,x_i)-\mu(x_s,Ax_i),x_j)
\overline{\langle\mu(x_s,x_i),x_j\rangle}  \\
&=\sum_{ijs}\left(\langle \mu(x_s,x_i),Ax_j\rangle -\langle \mu(Ax_s,x_i),x_j\rangle
-\langle \mu(x_s,Ax_i),x_j\rangle\right)\overline{\langle\mu(x_s,x_i),x_j\rangle}  \\
&=\sum_{ijrs} \left(\overline{\langle Ax_j,x_r\rangle}\langle\mu(x_s,x_i),x_r\rangle
-\langle Ax_s,x_r\rangle\langle\mu(x_r,x_i),x_j\rangle
-\langle Ax_i,x_r\rangle\langle\mu(x_s,x_r),x_j\rangle\right)\\
&\qquad\qquad \times \overline{\langle\mu(x_s,x_i),x_j\rangle}  \\
&=\sum_{ijrs}\overline{\langle\mu(x_i,x_j),x_s\rangle}
\langle\mu(x_i,x_j),x_r\rangle
\overline{\langle Ax_s,x_r\rangle}
-2\overline{\langle\mu(x_s,x_i),x_j\rangle}\langle\mu(x_r,x_i),x_j\rangle\langle Ax_s,x_r\rangle
\\
&=-2\sum_{ijrs}\overline{\langle L^\mu_{x_s}x_i,x_j\rangle}\langle L_{x_r}^\mu x_i,x_j\rangle
\langle Ax_s,x_r\rangle +\sum_{ijrs}\underbrace{\overline{\langle L_{x_i}^\mu x_j,x_s\rangle}}_{\overline{\langle x_j,L_{x_i}^{\mu*}x_s\rangle}} 
\underbrace{\langle L_{x_i}^\mu x_j,x_r\rangle}_{=\langle x_j,L^{\mu*}_{x_i}x_r\rangle} \overline{\langle Ax_s,x_r\rangle} \\
&=-2\sum_{irs} \langle L^\mu_{x_r} x_i,L^\mu_{x_s}x_i\rangle \langle x_s,Ax_r\rangle
+\sum_{irs}\langle L_{x_i}^{\mu*}x_s,L^{\mu*}_{x_i}x_r\rangle\overline{\langle Ax_s,x_r\rangle}\\
&=-2\sum_{irs} \langle L^\mu_{x_i} x_r,L^\mu_{x_i}x_s\rangle\overline{\langle Ax_r,x_s\rangle}
+\sum_{irs}\langle L_{x_i}^{\mu*}x_s,L^{\mu*}_{x_i}x_r\rangle\overline{\langle Ax_s,x_r\rangle}\\
&=-2\sum_i(L^{\mu*}_{x_i}L^\mu_{x_i},A)+\sum_i(L^{\mu}_{x_i}L^{\mu*}_{x_i},A),
\end{align*}
as wished. \EPf

\begin{cor}\label{cor:moment-jordan}
For all $x$, $y\in\C^n$ we have 
\[ \langle M_\mu x,y\rangle = -2\sum_{ij}\langle L_xx_i,x_j\rangle\overline{\langle L_yx_i,x_j\rangle}
+\sum_{ij}\overline{\langle L_{x_i}x_j,x\rangle}\langle L_{x_i}x_j,y\rangle. \]
\end{cor}

\subsection{The structure of solitons}

Next, we prove Proposition~\ref{prop:struct}. 
Let $\Ja=(\C^n,\mu)$ be a commutative algebra, where $\mu\in V_n$.  
We have seen in~(\ref{self-similar}) that $\mu$ is a soliton 
if and only if $D_\mu:=M_\mu-c_\mu I$ kills $\mu$, for 
$c_\mu=-||M_\mu||^2/||\mu||^2\in\R$ (recall that $I\cdot\mu=-\mu$));
owing to~(\ref{A-action}), this is equivalent to $D_\mu$ being a 
derivation of~$\mu$. Since $M_\mu^*=M_\mu$ and $c_\mu\in\R$, also $D_\mu^*=D_\mu$.  
For part~(a), it remains only to show that if $\mu$ is a 
soliton then some positive multiple of $D_\mu$ has rational eigenvalues. 
 
Without loss of generality, we may assume $D_\mu\neq0$.
Let $x_1,\ldots,x_n$ be an orthonormal basis of eigenvectors of $D_\mu$;
then it is also a basis of eigenvectors of $M_\mu$. 
Consider the 'structure constants' $\mu_{ij}^k$ given by 
$\mu(x_i,x_j)=\sum_k\mu_{ij}^kx_k$ for all $i$, $j$. 
Let $\Lh$ be the subspace of $\Lg$ consisting of endomorphisms
of $\C^n$ that are diagonal on that basis,  
let $\alpha_{ij}^k\in\Lh$ have matrix
$-E_{ii}-E_{jj}+E_{kk}$ in that basis, 
where $E_{ab}$ has a $1$ in the $(a,b)$-entry and $0$ elsewhere, 
and consider the subspace $F$ of $\Lh$ spanned by all  
$\alpha_{ij}^k$ with $\mu_{ij}^k\neq0$. 

Note that $A\in\Lh$ is a derivation of $\mu$  
if and only if $\mu(x_i,x_j)$ lies in the $(a_i+a_j)$-eigenspace
of~$A$, for all~$i$, $j$;
this is equivalent to $a_k=a_i+a_j$ whenever $\mu_{ij}^k\neq0$.
Therefore $\mathrm{Der}(\mu)\cap\Lh=F^\perp\cap\Lh$. 
Owing to $M_\mu\in \Lh$ and Lemma~\ref{lem:moment}(b), 
we have $M_\mu\in F$.
Now let $P:\Lh\to F$ be orthogonal projection. Applying~$P$
throughout the equation $M_\mu=c_\mu I+D_\mu$ yields
\[ M_\mu = c_\mu P(I). \]
Therefore
\[ -\frac1{c_\mu}D_\mu=I-\frac1{c_\mu}M_\mu=I-P(I). \]
Since $F$ is spanned by matrices with integer coefficients, $P(I)$
has rational coefficients, and this finishes the proof of part~(a).

We next address part~(b) of Proposition~\ref{prop:struct}.
Due to~\cite[p.~869]{jacobson1949}, the radical $\Jn$ 
of $\Ja$ is a characteristic ideal of $\Ja$, thus $D_\mu$-invariant.
Next, using results of Mostow and Auslander-Brezin,
we show there is a maximal semisimple subalgebra $\Js$ of $\Ja$
which is $D_\mu$-invariant. First we note that $D_\mu$ is semisimple,
since it is Hermitean. Therefore the one dimensional Lie algebra $\Ld$ 
of derivations of $\Ja$ generated by $D_\mu$ is completely reducible.
Its algebraic hull $\Ld^\#$ consists of derivations~\cite[(1.2)]{AB},
and it is also completely reducible~\cite[(1.5)]{AB}. Now the associated
connected algebraic subgroup of $GL(\Ja)$ is a completely reducible group
of automorphisms of~$\Ja$~\cite[(1.6)]{AB}. But a completely
reducible group of automorphisms of a Jordan algebra preserves some maximal
semisimple subalgebra~$\Js$ \cite[p.~215]{M}, and hence the same is true 
of $D_\mu$. 

Finally, we prove that $D_\mu|\Js=0$. Since $D_\mu$ is a Hermitean
endomorphism of $\C^n$, it is semisimple and has real eigenvalues.
Let $x\in\C^n$ be an eigenvector of $D_\mu$ with associated eigenvalue 
$d\in\R$. The derivation property implies that 
$x^k$ is an eigenvector of $D_\mu$ 
with eigenvalue $kd$ for every integer $k\geq2$.
By finite-dimensionality, if $d\neq0$ then $x$ is nilpotent;
it thus follows that $x\in\Jn$. 
This shows that the only eigenvalues of $D_\mu$  
that can occur for eigenvectors in $\Js$ are zero, and finishes the 
proof of the proposition. 
 
 \section{Semisimple Jordan algebras}

 In this section  we show that the
semisimple Jordan algebras realize the minimal value of the 
energy $E_n$ for all $n$, up to a change of basis
(Proposition~\ref{min-fn});
Theorem~\ref{ss-distinguished} is a consequence.
A good reference for Jordan-related results in this section 
is~\cite[Kap.~VIII]{braun-koecher}.

   \subsection{Peirce decomposition}

    Let $\mathfrak A$ be a complex semisimple Jordan algebra.
    Then there exists a \emph{complete orthogonal system of
      idempotents}, or \emph{Jordan frame}, that is, a set
    $\{e_1,\ldots,e_r\}$ of idempotents of $\mathfrak A$
    such that $e_1+\cdots+e_r$ is the identity element of $\mathfrak A$,
    and $e_ie_j=0$ if $i\neq j$.
    A Jordan frame is unique up to 
    an automorphism of $\mathfrak A$, and the number $r\geq2$
of elements in a Jordan frame 
    is called
    the \emph{degree} of $\mathfrak A$. 

A Jordan frame as above 
    gives rise to a canonical decomposition of $\mathfrak A$
into a vector space direct sum, as follows.
    The eigenvalues of an idempotent $e_k$ can only be $0$, $\frac12$ and $1$,
    and we denote the corresponding eigenspaces by
    $\mathfrak A_0(e_k)$, $\mathfrak A_{\frac12}(e_k)$
and $\mathfrak A_1(e_k)$, respectively.  The \emph{Peirce 
decomposition} of $\mathfrak A$ is
    \[ \mathfrak A=\oplus_{i\leq j}\mathfrak A_{ij}, \]
    where
    \[ \mathfrak A_{ii}=\mathfrak A_1(e_i) \]
    for all~$i$ and
\[   \mathfrak A_{ij}=\mathfrak A_{\frac12}(e_i)
  \cap \mathfrak A_{\frac12}(e_j) \]
for $i\neq j$. 

\subsection{Killing form}

Next, assume in addition to the above that $\mathfrak A$ is simple algebra
of dimension $n$.
Fix a Jordan frame and the associated Peirce decomposition.
It is known that if $r\geq3$, then 
$\dim\mathfrak A_{ij}=d\in\{1,2,4,8\}$ for all $i$, $j$.
Also, in case $r=2$ we put $d=n-2$. Now we have
\begin{equation}\label{n}
  n=r+\frac{r(r-1)}{2}d.
  \end{equation}

  An element $a\in\mathfrak A$ is called \emph{regular}
if $\C[a]:=\mathrm{span}_{\mathbb C}\{1,a,\ldots,a^{r-1}\}$ is
$r$-dimensional (equivalently, $\dim\C[a]$ is maximal
among $\dim\C[b]$ for all $b\in\mathfrak A$).
The \emph{reduced trace} of a regular element
$a\in \mathfrak A$ is
  \[ \mathrm{tr}(x)=\mathrm{Tr}(L|_{\mathbb C[x]}), \]
  and the function ``reduced trace'' can be uniquely extended
  to a linear map $\mathrm{tr}:\mathfrak A\to\C$ 
(compare~\cite[chap.2, \S2]{FK}). 
Since $\mathfrak A$ is
  simple, we have
  \begin{equation}\label{tr-Tr}
    \mathrm{tr}(ab)=\frac rn\tau(a,b).
    \end{equation}

The \emph{Killing form} of $\mathfrak A$ is the symmetric bilinear form
$K:\mathfrak A\times\mathfrak A\to\C$ defined by 
\[ K(a,b) = \mathrm{Tr}(L_aL_b), \]
for $a$, $b\in \mathfrak A$. It follows from the trace formula
that~\cite[Satz~9.4, Kap.~VIII]{braun-koecher}
\begin{equation}\label{killing}
  K(a,b)=\left(1+(r-2)\frac d4\right)\mathrm{tr}(ab)+\frac d4\mathrm{tr}(a)\mathrm{tr}(b).
  \end{equation}

\subsection{Moment map}

Suppose $\mathfrak A$ is a complex semisimple Jordan algebra
given by $\mu\in\mathcal J_n$. Then the trace form 
$\tau$ is a nondegenerate symmetric
bilinear form. Also, there exists an Euclidean real form
$\mathfrak A_0$ of $\Ja$, that is, $\mathfrak A_0$ is
a formally real Jordan algebra such that
$\mathfrak A_0\otimes_{\mathbb R}\C=\mathfrak A$~(\cite[Satz~5.6, Kap.~IX]{braun-koecher} or~\cite[Thm.~8.5.2]{FK}). 
The restriction of $\tau$
to $\mathfrak A_0$ is a positive-definite Euclidean inner
product~(\cite[Satz~3.4, Kap.~IX]{braun-koecher} or~\cite[Prop.~3.1.5]{FK});
we extend it to an Hermitean inner product on $\mathfrak A$,
denoted by $\langle,\rangle$.

\begin{lem}\label{adj-left-mult}
  Let $a\in\mathfrak A$. Then the adjoint $L_a^*$
  of $L_a:\mathfrak A\to\mathfrak A$
  with respect to $\langle,\rangle$ is $L_{\bar a}$, where
  $\bar a\in\mathfrak A$ is the complex conjugate of~$a$ over $\mathfrak A_0$.
  In particular, $L_a$ is a Hermitean operator for $a\in\mathfrak A_0$. 
\end{lem}

\Pf We compute for $b$, $c\in\mathfrak A$ that
\[ \langle L_ab,c\rangle = \tau(ab,\bar c)=\tau(b,a\bar c)
  =\tau(b,\overline{L_{\bar a}c})=\langle b,L_{\bar a}c\rangle, \]
since $L_a$ is self-adjoint for~$\tau$, 
which proves the statement. \EPf

\medskip

The next result shows that the moment matrix of a
complex semisimple Jordan algebra $\mathfrak A$
is essentially given by its Killing form.

\begin{prop}\label{moment-killing}
  Let $M_\mu$ be the Hermitean matrix which is the moment matrix of $\mathfrak A=(\C^n,\mu)$, and let $K_\mu$ be the Killing form of $\mathfrak A$. Then
  \[ \langle M_\mu a,b\rangle = -K_\mu(a,\bar b) \]
  for all $a$, $b\in \mathfrak A$.
\end{prop}

\Pf Let $\{x_1,\ldots,x_n\}$ be a orthonormal basis of $\C^n$ with respect to
$\langle,\rangle$ which is contained in $\mathfrak A_0$. 
Using Proposition~\ref{moment-jordan} and Lemma~\ref{adj-left-mult},
we can write
\begin{align*}
  \langle M_\mu a,b\rangle &= -\sum_i\langle (L_{x_i}^\mu)^2a,b\rangle\\
                           &=-\sum_i\langle L_{x_i}^\mu a,L_{x_i}^\mu b\rangle\\
   &=-\sum_i\langle L_a^\mu{x_i},L_b^\mu{x_i} \rangle\\
                           &=-\sum_i\langle L_{\bar b}^\mu L_a^\mu{x_i},{x_i}\rangle\\
                           &=-\mathrm{Tr}(L^\mu_{\bar b}L^\mu_a)\\
                           &=-K_\mu(a,\bar b),
\end{align*}
as desired. \EPf

\begin{lem}\label{moment-matrix-peirce-basis}
  Let $\mathfrak A$ be a complex simple Jordan algebra
  represented by $\mu\in\mathcal J_n$. Denote by $r$ and $d$, respectively,
  the
  degree of $\mathfrak A$ and the dimension of the off-diagonal
  Peirce components $\mathfrak A_{ij}$ ($i\neq j$), where we
  fix a Jordan frame $(e_1,\ldots,e_r)$ contained in an Euclidean
  real form. Fix also a Hermitian
  product on $\C^n$ by putting (see~(\ref{tr-Tr}))
  \[ \langle x,y\rangle :=\frac rn\mathrm{Tr}(L^\mu_{x\bar y})=\mathrm{tr}^\mu(x\bar y) \]
  for $x,y\in\C^n$. 
Then the Jordan
frame is a orthonormal set and can be completed to a orthonormal
basis of $\C^n$ by adding elements $e_{ij}^k\in\mathfrak A_{ij}$,
$k=1,\ldots,d$, for all  $i\neq j$.
In this basis, the moment matrix of $\mathfrak A=(\C^n,\mu)$ is given by
  \[ M_\mu = -\alpha I -\frac d4 N \]
  where
  $\alpha=1+(r-2)\frac d4$ and 
  \[
N=\underset{\hspace{.3in}r\hspace{.4in}n-r}{\left(\begin{array}{c|c}
  \begin{matrix}
    1 &\cdots& 1 \\
    \vdots &&\vdots\\
  1 &\cdots &  1
  \end{matrix}
  &  \bigzero\\
\hline\\
    \bigzero & \bigzero 
\end{array}\right).}
\]
\end{lem}

\Pf Consider the Peirce decomposition
$\mathfrak A=\oplus_{i\leq j}\mathfrak A_{ij}$
with respect to a Jordan frame $(e_1,\ldots,e_r)$.
Then
\[ \mathfrak A_0(e_k)=\oplus_{i,j\neq k}\mathfrak A_{ij},\
  \mathfrak A_{\frac12}(e_k)=\oplus_{i\neq k}\mathfrak A_{ik},\ \mbox{and}\
  \mathfrak A_1(e_k)=\mathfrak A_{kk}, \]
where $\dim\mathfrak A_{ij}=d$ for all $i\neq j$, and
$\mathfrak A_{kk}=\C\cdot e_k$ for all $k$. Using $e_k^2=e_k$,
(\ref{n}) and~(\ref{tr-Tr}),
we see immediately that
\[ \langle e_k,e_k\rangle=\mathrm{tr}^\mu(e_k)=\frac rn \mathrm{Tr}(L^\mu_{e_k}) =
  \frac rn \left(\frac12 (r-1)d + 1\right) = 1 \]
for all $k$. Moreoever, using
\begin{gather*}
  \mathfrak A_{ij}\cdot(\mathfrak A_{ii}+\mathfrak A_{jj})\subset
  \mathfrak A_{ij},\ \mathfrak A_{ij}\cdot\mathfrak A_{ij}\subset
  \mathfrak A_{ii}+\mathfrak A_{jj}\\
  \mathfrak A_{ij}\cdot\mathfrak A_{jk}\subset\mathfrak A_{ik},\
  \mathfrak A_{ij}\cdot\mathfrak A_{ik}\subset\mathfrak A_{jk},
\end{gather*}
for mutually different $i$, $j$, $k$, we see that, for 
$e_{ij}^k\in\mathfrak A_{ij}$ ($i\neq j$),
\[ \mathrm{tr}^\mu(e_{ij}^k)=\frac rn\mathrm{Tr}(L^\mu_{e_{ij}^k})=0. \]
Finally, formula~(\ref{killing}) and Proposition~\ref{moment-killing}
yield:
\[ \langle M_\mu x,y\rangle = -\alpha \langle x,y\rangle -\frac d4
  \mathrm{tr}^\mu(x)\mathrm{tr}^\mu(\bar y). \]
Therefore
\begin{gather*}
  \langle M_\mu e_i,e_j\rangle = -\alpha\delta_{ij}-\frac d4, \\
  \langle M_\mu e_\ell,e_{ij}^k\rangle = 0,
\end{gather*}
and
\[ \langle M_\mu e_{ij}^k,e_{rs}^t\rangle = -\alpha \delta_{ir}\delta_{js}
  \delta_{kt}, \]
as wished. \EPf
  
\begin{rmk}
 Using the calculation in the proof of Lemma~\ref{moment-matrix-peirce-basis}, it is easy to see that 
\[ E_n(\mu)=\frac{\mathrm{Tr} M_\mu^2}{(\mathrm{Tr}M_\mu)^2}
=\frac{n\alpha^2+\frac{r^2d^2}{16}+\frac{r\alpha d}2}{(n\alpha+\frac{rd}4)^2}
=\frac{n\alpha^2+\frac{r\alpha d}2+\frac{r^2d^2}{16}}{n\alpha^2+\frac{r\alpha d}2+\frac{r^2d^2}{16n}}\cdot\frac 1n> \frac 1n, \]
where $n\geq2$. 
In Proposition~\ref{min-fn}, we will see that the Peirce basis 
can be slightly changed to lower the value of $E_n$.
\end{rmk}

Theorem~\ref{ss-distinguished} is derived from the following
result. 

  \begin{prop}\label{min-fn}
    For all $n$, the minimum value of $E_n:\mathcal J_n\to\R$
    is $1/n$, and this value is attained only at semisimple
    Jordan algebras. Conversely, if $\mu$ is semisimple then
    $E_n(g\cdot\mu)=1/n$ for some $g\in G$. 
  \end{prop}

  \Pf Note that $E_n(\mu)=\frac1n$ if and only if $M_\mu=c_\mu I$.
Indeed if $c_1,\ldots, c_n$ are the (real) eigenvalues of $M_\mu$,
and we want to minimize the value of $c_1^2+\cdots+c_n^2$ subject to
the condition that $c_1+\cdots+c_n=1$, we immediately obtain that 
$c_1=\cdots=c_n=1/n$. 
  
We first prove that the moment matrix of a complex
  semisimple Jordan algebras is scalar, up to a change of basis. 
  It is enough to consider a complex simple Jordan algebra $\mu$.  
We start with the basis $\{e_i\}\cup\{e_{ij}^k\}$ of
the Euclidean real form as in Lemma~\ref{moment-matrix-peirce-basis}
and the associated moment matrix $M_\mu=-\alpha I -\frac d4 N$.
Let $\{f_1,f_2,\ldots,f_r\}$ be a (positively oriented) real orthonormal basis
of $\mathrm{span}_{\mathbb R}(e_1\ldots,e_r)$ such that its
first element is parallel to the identity element $1=e_1+\ldots+e_r$.
Then $\{f_i\}\cup\{e_{ij}^k\}$ is an orthonormal basis of $\R^n$.
Let $k\in SO(n)$, $k=\left(\begin{smallmatrix}k'&0\\0&I\end{smallmatrix}\right)$,
where $k'\in SO(r)$, be such that $k$ maps $f_i\mapsto e_i$, $e_{ij}^k\mapsto e_{ij}^k$,
and put $\nu=k\cdot\mu$. Since the moment map is $U(n)$-equivariant,
we have
\begin{equation}\label{mnu}
 M_\nu=kM_\mu k^{-1}=-\alpha I-\frac d4 \left(\begin{array}{c|c}
  \begin{matrix}r&&&\\&0&&\\&&\ddots&\\&&&0\end{matrix}&\bigzero\\ \hline\\\bigzero&\bigzero
  \end{array}\right). 
\end{equation}
We next show that it is possible to rescale the basis (in fact, only $f_1$)
so that the moment matrix becomes scalar. 

We need to make some remarks about the $\nu_{ij}^k$.
Denote by $\{f_1,\ldots,f_n\}$ the orthonormal basis of $\R^n$ 
above constructed so that $f_1$ is parallel to $1$,
$\{f_1,\ldots, f_r\}$ spans $\mathfrak A_{11}+\ldots+\mathfrak{A}_{rr}$
and $\{f_{r+1},\ldots,f_n\}$ spans $\oplus_{ij}\mathfrak A_{ij}$.  
Let $x\in \R^n$ be a unit vector. Then ($x^2=\mu(x,x)$)
\[ x^2 = \frac{\langle x^2,1\rangle}{||1||^2}1+ y,\quad y\perp x. \]
Note that $\langle x^2,1\rangle=||x||^2=1$ and $||1||^2=\frac rn\mathrm{Tr}(I)=r$. It follows that
\begin{equation}\label{x^2}
  x^2 = \frac 1{\sqrt r} f_1 +y,\quad y\perp f_1.
\end{equation}
We also have $\langle \mu(f_i,f_j),1\rangle=\langle f_i,f_j\rangle =0$
if $i\neq j$. Since $\nu(e_i,e_j)=k(\mu(f_i,f_j))$, 
it follows that
\begin{equation}\label{form-nu}
  \nu = \frac1{\sqrt r}\sum_{i=1}^ne_1'e_i'\otimes e_i
  +\frac1{\sqrt r}\sum_{i=2}^n e_i^{\prime 2}\otimes e_1 +\mbox{terms
    not involving $e_1$, $e_1'$;}
  \end{equation}
here we denote by $e_1',\ldots,e_n'$ the dual basis of $e_1,\ldots,e_n$. 
  
Next we consider the following one-parameter deformation of $\nu$, and show
that the moment matrix is scalar for some value of the parameter.
Set $\nu_t=g_t^{-1}\cdot\nu$, where
\[ g_t=\exp(-tE_{11})=\begin{pmatrix}e^{-t}&&&\\
    &1&&\\ &&\ddots& \\ &&& 1\end{pmatrix}. \]
We can write $\nu=\sum_{ijk}\nu_{ij}^kv^{ij}_k$, where
$v^{ij}_k=e_i'e_j'\otimes e_k$ and 
the sum runs through $i\leq j$, and then
\[ ||v_k^{ij}||^2=2\ \mbox{($i\neq j$) and}\  ||v_k^{ii}||^2=1. \]
Also, $v^{ij}_k$ is a weight vector of $S^2(\C^{n*})\otimes\C^n$
of weight $\alpha_{ij}^k=-\theta_i-\theta_j+\theta_k$, with respect
to the Cartan subalgebra of $\Lg$ consisting of diagonal matrices, 
where we denote by $\theta_i$ the projection onto the $i$th-diagonal
entry. Now
\[ \nu_t=\sum_{ijk}e^{t\alpha_{ij}^k(E_{11})}\nu_{ij}^kv^{ij}_k \]
and
\[ E_{\ell\ell}\cdot\nu_t=\sum_{ijk}e^{t\alpha_{ij}^k(E_{11})}\alpha_{ij}^k(E_{\ell\ell})\nu_{ij}^kv_k^{ij}. \]
It follows that, for $\ell\geq2$, 
\begin{align}\nonumber
  (M_{\nu_t},E_{\ell\ell})&=\langle E_{\ell\ell}\cdot\nu_t,\nu_t\rangle \\ \nonumber
                          &= \sum_{ijk}e^{2t\alpha_{ij}^k(E_{11})}\alpha_{ij}^k(E_{\ell\ell})|\nu_{ij}^k|^2||v^{ij}_k||^2\\ \nonumber
                          &= 0+\frac 1r e^{2t}(-2)\cdot 1 +\underbrace{\sum_{ijk\neq1}\alpha_{ij}^k(E_{\ell\ell})|\nu_{ij}^k|^2||v^{ij}_k||^2}_{=:c_\ell}\\ \label{mnutell}
                          &=-\frac 2re^{2t}+c_\ell,
  \end{align}
where we have used the form~(\ref{form-nu}). 
Since $(M_{\nu_0},E_{\ell\ell})=-\frac 2r+c_\ell$ is independent of 
$\ell\geq2$, owing to~(\ref{mnu}) 
we deduce that $c_2=\cdots=c_n=c$ for some $c\in\R$. 

We next describe the first entry of $M_{\nu_t}$, again by
using~(\ref{form-nu}):
\begin{align}\nonumber
  (M_{\nu_t},E_{11}) & =  \sum_{ijk}e^{2t\alpha_{ij}^k(E_{11})}\alpha_{ij}^k(E_{11})|\nu_{ij}^k|^2||v^{ij}_k||^2\\ \nonumber
                     &=\frac 1r e^{-2t}(-1)(1+(n-1)2)+\frac 1r e^{2t}\cdot 1\cdot ((n-1)\cdot1)
  \\ \label{mnute11}
                     &= -\frac{2n-1}r e^{-2t}+\frac{n-1}r e^{2t}.
                       \end{align}
                       Use~(\ref{mnutell}) and~(\ref{mnute11}) to
                       investigate ($\ell\geq2$)
                       \[ \gamma(t)=(M_{\nu_t},E_{11})-(M_{\nu_t},E_{\ell\ell})=\frac{n+1}re^{2t}-
                         \frac{2n-1}re^{-2t}-c. \]
Since $\lim_{t\to-\infty}\gamma(t)=-\infty$ and
$\lim_{t\to+\infty}\gamma(t)=+\infty$, there is $t_0\in \R$ such that 
$\gamma(t_0)=0$. We have shown that all diagonal 
entries of $M_{\nu_{t_0}}$ are equal. The last step is to show that 
the off-diagonal entries of $M_{\nu_{t_0}}$ (indeed of $M_{\nu_t}$ for all 
$t$) vanish.                         

Since $M_\nu$ is a diagonal matrix, it suffices to show that ($1\leq p<q\leq n$):
\begin{equation}\label{off-diag-mnut}
 \frac d{dt}(M_{\nu_t},E_{pq}+E_{qp})=2\Re\langle (E_{pq}+E_{qp})\cdot\nu_t,E_{11}\cdot\nu_t\rangle 
\end{equation}
vanishes. 

Note that 
\[ E_{11}\cdot\nu_t=-\frac{e^{-t}}{\sqrt r}\sum_{i=1}^nv^{1i}_i
+\frac{e^t}{\sqrt r}\sum_{i=2} ^nv^{ii}_1 \]
so we need only know the $v^{1i}_i$- and $v^{ii}_1$-components
of~$(E_{pq}+E_{qp})\cdot\nu_t$. 

We compute
\begin{gather*}\label{1} \nonumber
\langle((E_{pq}+E_{qp})\cdot\nu_t)(e_1,e_i),e_i\rangle=-\delta_{p1}\langle
\nu_t(e_q,e_i),e_i\rangle,\\
\langle(E_{11}\cdot\nu_t)(e_1,e_i),e_i\rangle =-\langle\nu_t(e_1,e_i),e_i\rangle,
\end{gather*}
and, for $i=2,\ldots,n$, 
\begin{gather*}\label{2}
\langle((E_{pq}+E_{qp})\cdot\nu_t)(e_i,e_i),e_1\rangle=\delta_{p1}\langle
\nu_t(e_i,e_i),e_q\rangle,\\
\langle(E_{11}\cdot\nu_t)(e_i,e_i),e_1\rangle =\langle\nu_t(e_i,e_i),e_1\rangle.
\end{gather*}
Plugging these formulae into~(\ref{off-diag-mnut}) already yields zero,
unless $p=1$, for which
\begin{equation}\label{off-diag-mnut2}
 \frac12\frac d{dt}(M_{\nu_t},E_{1q}+E_{q1})=
\Re\sum_{i=2}^n\langle\nu_t(e_q,e_i),e_i\rangle\overline{\langle\nu_t(e_1,e_i),e_i\rangle}+\langle\nu_t(e_i,e_i),e_q\rangle\overline{\langle\nu_t(e_i,e_i),e_1\rangle}. 
\end{equation}
Note that, using~$i$, $q\geq2$,  
\begin{gather*}
\langle\nu_t(e_i,e_i),e_q\rangle = \langle \nu(e_i,e_i),e_q\rangle,\\
\langle\nu_t(e_i,e_i),e_1\rangle = e^t\langle \nu(e_i,e_i),e_1\rangle,\\
\langle\nu_t(e_q,e_i),e_i\rangle = \langle \nu(e_i,e_i),e_q\rangle,\\
\langle\nu_t(e_1,e_i),e_i\rangle = e^{-t}\langle \nu(e_i,e_i),e_1\rangle,\\
\end{gather*}
so 
\begin{align*}
 \frac d{dt}(M_{\nu_t},E_{1q}+E_{q1})&=
\sum_{i=2} ^n e^{-t}\langle \nu(e_i,e_i),e_q\rangle
\langle \nu(e_i,e_i),e_1\rangle+
e^t\langle \nu(e_i,e_i),e_q\rangle
\langle \nu(e_i,e_i),e_1\rangle\\
&=\sum_{i=2} ^n\frac1{\sqrt r}(e^t+ e^{-t})\langle \nu(e_i,e_i),e_q\rangle\\
&=\frac1{\sqrt r}(e^t+e^{-t})\left\langle\sum_{i=2} ^n\nu(e_i,e_i),e_q\right\rangle. 
\end{align*}
Finally, we claim that the basis $f_1,\ldots,f_n$ can be chosen 
so that $\sum_{i=2}^n\mu(f_i,f_i)\in\R\cdot 1$. This 
will imply $\sum_{i=2} ^n\nu(e_i,e_i)\in\R\cdot e_1\perp e_q$,
and hence $\frac d{dt}(M_{\nu_t},E_{1q}+E_{q1})=0$. 

The claim is proved in two steps. In the first step, recall the 
idempotents $e_1\ldots,e_r$ with $e_ie_j=0$ for $i\neq j$ 
and $e_1+\ldots+e_r=1$. Set $f_1=\tfrac1{\sqrt r}$ and 
\[ f_i=\frac1{\sqrt{i(i-1)}}(e_1+\ldots+e_{i-1}-(i-1)e_i) \]
for $i=2,\ldots, r$. Then a quick calculation yields
\[ \sum_{i=2} ^r f_i^2 =\left(1-\frac 1r\right)\cdot 1. \]

In the second step, we invoke~\cite[Satz~9.1, Kap.~VIII]{braun-koecher}.
It says that $u^2=\frac12(e_i+e_j)$ for all $u\in\mathfrak A_{ij}$ ($i\neq j$) with $||u||=1$. Finally, 
\[ \sum_{i=r+1}^n f_i^2=\sum_{1\leq i<j\leq r} \frac d2 (e_i+e_j)
=\frac{(r-1)d}2\sum_{i=1}^r e_i = \frac{(r-1)d}2\cdot 1, \]
and hence 
\[ \sum_{i=1}^n f_i^2=(r-1)\left(\frac 1r + \frac d2\right)\cdot 1, \]
as wished. 

In the remainder of the proof, we show that, conversely, 
if $M_\mu$ is a scalar matrix then $\mathfrak A=(\C^n,\mu)$ is 
semisimple. 

Write $\mathfrak A=\mathfrak N+\mathfrak S$ (direct sum of vector 
spaces), where $\mathfrak N$ is the 
radical of $\mathfrak A$ and $\mathfrak S$ is a semisimple 
Jordan algebra. Suppose, by contradiction, that $\mathfrak N\neq0$. 
Then 
\[ \mathfrak N^{[0]}=\mathfrak N,\ \mathfrak N^{[k+1]}=(\mathfrak N^{[k]})^3,
\] 
for $k\geq0$ defines a decreasing sequence of ideals of 
$\mathfrak A$~(\cite[Lemma~2.2]{penico}; see
also~\cite[Lem.~3, \S3, ch.4]{ZSSS}). Nilpotency of $\mathfrak N$ 
yields a minimal $k_0\geq0$ such that $\mathfrak N^{[k_0+1]}=0$. Let 
$\mathfrak B=\mathfrak N^{[k_0]}\neq0$. Then  
\[ \mathfrak B_0=\mathfrak B,\ \mathfrak B_{k+1}=\mathfrak A\mathfrak B_k^2+\mathfrak B_k^2, \]
for $k\geq0$, is a decreasing sequence of ideals 
of~$\mathfrak A$~\cite[Lemma~2.2]{penico}. By~\cite[Theorem~2.5]{penico},
there is a minimal $k_1\geq1$ such that $\mathfrak B_{k_1}\subset\mathfrak B^2$. 
Now there are two cases.

If $\mathfrak B_{k_1}\neq0$, then $\mathfrak C=\mathfrak B_{k_1}$ is a 
nonzero ideal of $\mathfrak A$ with
\[ \mathfrak C^2=\mathfrak B_{k_1}^2\subset\mathfrak B^2\mathfrak B=\mathfrak B^3
=(\mathfrak N^{[k_0]})^3=\mathfrak N^{[k_0+1]}=0. \]

If  $\mathfrak B_{k_1}=0$, then $\mathfrak C=\mathfrak B_{k_1-1}$ is a 
nonzero ideal of $\mathfrak A$ with
\[ \mathfrak C^2=\mathfrak B_{k_1-1}^2=0 \]
since $0=\mathfrak B_{k_1}=\mathfrak A\mathfrak B_{k_1-1}^2+\mathfrak B_{k_1-1}^2$. 

In any case, we have found a nonzero ideal $\mathfrak C$ of $\mathfrak A$ 
with $\mathfrak C^2=0$. We can now finish the proof. 
Let $x_1,\ldots,x_m$ be a orthonormal basis of $\mathfrak C$, 
and let $y_1,\ldots,y_{n-m}$ be a orthonormal basis of $\mathfrak C^\perp$. 
For any $x\in\mathfrak C$, owing to Corollary~\ref{cor:moment-jordan}
and the facts that $L_x^\mu x_i=0$ for all $i$ and
that the image of $L_x^\mu$ is contained in $\mathfrak C$, 
we have:
\begin{align*}
0&>\langle M_\mu x,x\rangle\\
&=-2\sum_{ij}|\langle L_x^\mu y_i,x_j\rangle|^2
+\sum_{ij}|\langle L_{x_i}^\mu y_j,x\rangle|^2
+\sum_{ij}|\langle \underbrace{L_{y_i}^\mu x_j}_{=L^\mu_{x_j}y_i},x\rangle|^2
+\sum_{ij}|\langle L_{y_i}^\mu y_j,x\rangle|^2.
\end{align*}
We make $x=x_k$ and sum over $k=1,\ldots,m$ to obtain
\begin{align*}
0&>\sum_k\langle M_\mu x_k,x_k\rangle\\
&=-2\sum_{ijk}|\langle L_{x_k}^\mu y_i,x_j\rangle|^2
+2\sum_{ijk}|\langle L_{x_i}^\mu y_j,x_k\rangle|^2
+\sum_{ijk}|\langle L_{y_i}^\mu y_j,x_k\rangle|^2\\
&=\sum_{ijk}|\langle L_{y_i}^\mu y_j,x_k\rangle|^2\\
&\geq0,
\end{align*}
a contradiction. Hence $\mathfrak N=0$, as desired. \EPf

\section{The maximal value of $E_n$}\label{max}

Since the energy $E_n$ is constant along rays in $V_n\setminus\{0\}$, 
it attains a maximum value.
In this section, we determine those points of maxima and 
prove Theorem~\ref{nilp-comp}.

We introduce two important complex Jordan algebras. 
We give them names following the analogy with Lie algebras:
\begin{enumerate}
\item The \emph{Heisenberguian Jordan algebra} $\mu_{Heis}$ has a basis
$\{n_1,\ldots, n_n\}$ ($n\geq2$) satisfying $n_1^2=n_2$, and the 
other products equal to zero. This is a nilpotent 
Jordan algebra. 
\item The \emph{hyperbolic Jordan algebra} $\mu_{hyp}$ 
has a basis $\{e,n_1,\ldots,n_{n-1}\}$ ($n\geq2$) satisfying
$e^2=e$, $en_i=\frac12n_i$ for all $i$, and 
the other products equal to zero. This algebra has been 
considered in~\cite[Teorema~4.65]{martin},
where it was shown that it is rigid.
\end{enumerate}

\begin{prop}\label{deg}
Every Jordan algebra of dimension at least two which is not isomorphic 
to $\mu_{hyp}$ degenerates to~$\mu_{Heis}$. Further, the only 
Jordan algebras in $\mathcal J_n$ ($n\geq2$) for which the only non-trivial degeneration 
is to the trivial Jordan algebra (all products zero) are 
$\mu_{Heis}$ and~$\mu_{hyp}$.
\end{prop}

\Pf Let $\mu\in\mathcal J_n$. Suppose there is $x_1\in\C^n$ such that 
$x_1$, $x_2:=\mu(x_1,x_1)$ are linearly independent; complete 
this set to a basis $x_1,\ldots, x_n$. Define $g_t\in G$ 
by setting $g_tx_1=tx_1$, $g_tx_i=t^2x_i$ for $i=2,\ldots, n$,
and put $\mu_t=g_t^{-1}\cdot\mu$. Then $\mu_t(x_1,x_1)=x_2$ for all
$t\neq0$. Moreover,  
for $(i,j)\neq(1,1)$ there is $m\in\{3,4\}$ such that   
\[ \mu_t(x_i,x_j)=t^{m-1}\mu_{ij}^1x_1+t^{m-2}(\mu_{ij}^2x_2+\cdots
+\mu_{ij}^nx_n)\to 0 \]
as $t\to 0$. In other words, $\mu\to\mu_{Heis}$. 

In case there is no $x_1$ as a above, there is a non-zero linear map
$\ell:\C^n\to\C$ such that $\mu(x,x)=\ell(x)x$. 
By polarization, $\mu(x,y)=\frac12(\ell(x)y+\ell(y)x)$.
Let $n_1,\ldots,n_{n-1}$ be a basis of $\ker\ell$ and choose
$e$ such that $\ell(e)=1$. Then $\mu(e,e)=e$, $\mu(e,n_i)=\frac12n_i$
and $\mu(n_i,n_j)=0$ for all~$i$, $j$, which shows that 
$\mu\cong\mu_{hyp}$. 

 Further, $\dim\mathrm{Der}(\mu_{hyp})=n^2-n$
and $\dim\mathrm{Der}(\mu_{Heis})=n+(n-2)(n-1)=n^2-2n+2$; in fact,
if $d\in\mathrm{Der}(\mu_{Heis})$, then $d(n_1)$ is arbitrary, 
$d(n_i)\in\mathrm{span}(n_2,\ldots,n_n)$ for $i\geq2$, and
$d(n_2)=2n_1d(n_1)$. Since 
 $\dim\mathrm{Der}(\mu_{Heis})\leq\dim\mathrm{Der}(\mu_{hyp})$
for $n\geq2$, it follows that 
$\mu_{hyp}\not\to\mu_{Heis}$~\cite[p.~284]{kashuba-martin2014}. 
Also, $\mu_{Heis}\not\to\mu_{hyp}$ because 
$\dim\mu_{hyp}(\C^n,\C^n)=n>1= \dim\mu_{Heis}(\C^n,\C^n)$
(alternatively, $\mu_{Heis}$ is associative, but 
$\mu_{hyp}$ is not). \EPf

\medskip

Since $\mu_{hyp}$ is not a nilpotent Jordan algebra, 
it follows from Proposition~\ref{deg} that $\mu_{Heis}$ lies
in the closure of every $G$-orbit in $\mathcal N_n$,
which proves the first assertion of Theorem~\ref{nilp-comp}.
The second assertion is a consequence of the following result. 

\begin{cor}
  The Jordan algebras $\mu_{Heis}$ and $\mu_{hyp}$ are solitons.
  Further, the maximum value of $E_n$ is $5$,
    and it is attained exactly at the $G$-orbit of~$\mu_{Heis}$.
  \end{cor}

  \Pf Suppose $\mu$ is a point of maximum of $E_n$.
This implies that $\mu$ is a critical point of $E_n$ on $G\cdot\mu$,
and hence on~$V_n$.  
  By Kirwan-Ness theory (Theorem~\ref{kirwan-ness}),
$\mu$ is a point of minimum of
  $E_n$ on $G\cdot\mu$, which implies that $E_n$ is constant
  along $G\cdot\mu$. Now every point in $G\cdot\mu$ is a
  point of minimum of $E_n$ on $G\cdot\mu$, so the $G$- and
$K$-orbits through~$\mu$ agree up to scaling, that
is $G\cdot\mu=\C^\times\cdot K\cdot\mu$, again by Kirwan-Ness. 
  It follows that the only
  possible degeneration of $\mu$ is to the trivial algebra. By
  Theorem~\ref{deg}, $\mu$ is isomorphic to one of $\mu_{Heis}$ or~$\mu_{hyp}$.
  A simple calculation yields:
  \[ M_{\mu_{Heis}}=-5I+\begin{pmatrix}3&&&&\\ &6&&&\\ &&5&& \\ &&&\ddots&
      \\ &&&&5\end{pmatrix},\quad E_n(\mu_{Heis})=5, \]
  and 
  \[ M_{\mu_{hyp}}=-\left(\frac{n+1}2\right)I+\begin{pmatrix}0&&&\\ &\frac{n+1}2&&\\ &&\ddots& \\ &&&\frac{n+1}2\end{pmatrix},\quad E_n(\mu_{hyp})=1. \]
  This finishes the proof. \EPf

\section{Stratification}

In this section, we review the Kirwan-Ness
stratification of the null cone in the setting of commutative 
algebras. We follow the formulation and notation 
of~\cite{lauret2003,lauret-will} (for which we refer the reader, regarding
the missing proofs below) and note that 
the results are exactly the same as for the case of 
skew-symmetric algebras. 

Consider again the action of $G=GL(n,\C)$ on 
$V_n=S^2(\C^{n*})\otimes\C^n$. Since
$\lim_{t\to0}  g_t^{-1}\cdot \mu = \lim_{t\to0} t\mu =0$ for $g_t=tI$, 
every commutative algebra 
degenerates to the trivial algebra, and therefore the null cone $N=V_n$.
Let $\Lh$ be the subalgebra of diagonal matrices
of $i\Lk$, choose a positive Weyl chamber $\Lh^+$ in $\Lh$
and denote its closure by $\overline{\Lh^+}$. 
The critical set of $E$ is $K$-invariant and it decomposes
into a finite union of disjoint closed subsets
$\{C_\beta\}_{\beta\in\mathcal B}$, where $\mathcal B$ is a
finite subset of $\overline{\Lh^+}$, such that
$C_\beta$ is mapped under $m$ onto the adjoint
orbit $K\cdot\beta$ in~$i\Lk$. Let $S_\beta$ be the set of all points of 
$V_n\setminus\{0\}$ that flow into $C_\beta$ under
$-\nabla E_n$ (the stable manifold
of $C_\beta$). Then $S_\beta$ is $G$-invariant,
Zariski-locally closed, irreducible and non-singular, and we have:
\begin{equation}\label{strat}
 V_n\setminus\{0\}=\dot\bigcup_{\beta\in\mathcal B}S_\beta
\qquad\mbox{(disjoint union)}, 
\end{equation}
where
\begin{equation}\label{sb}
 \bar S_\beta \setminus S_\beta \subset \bigcup_{||\beta'||>||\beta||}S_{\beta'}. 
\end{equation}

Write $\mu=\sum_{ijk}\mu_{ij}^ke_i'e_j'\otimes e_k$,
for an orthonormal basis $e_1,\ldots,e_n$ of $\C^n$ 
and its dual basis $e_1',\ldots,e_n'$. 
For $\mu\in V_n\setminus\{0\}$, define
\[ \beta_\mu=\parbox[t]{5.5in}{the convex combination of smallest norm of the
    elements $\alpha_{ij}^k\in\Lh$ with $\mu_{ij}^k\neq0$.} \]
Recall that $\alpha_{ij}^k=-E_{ii}-E_{jj}+E_{kk}$, so $\mathrm{Tr}\beta_\mu=-1$.
Now another description is 
\[ S_\beta= \{\mu\in V_n\setminus\{0\}\;|\;\beta\ \mbox{is of maximal
    norm in $\{\beta_{g\cdot\mu}\;|\;g\in G\}$}\}, \]
and $\mathcal B=\{\beta\in\bar\Lh^+\;|\;S_\beta\neq\varnothing\}$. 
Define also 
\begin{equation*}\label{wb}
 W_\beta=\{\mu\in V_n\;|\;(\beta,\alpha_{ij}^k) \geq||\beta||^2\
  \mbox{for all $\mu_{ij}^k\neq0$} \},
\end{equation*}
that is, the sum of eigenspaces of $d\pi(\beta)$ with eigenvalues
$\geq||\beta||^2$, its subset
\[ Y_\beta=\{\mu \in W_\beta\;|\;(\beta,\alpha_{ij}^k)=||\beta||^2\
  \mbox{for at least one  $\mu_{ij}^k\neq0$} \}, \]
and
\[ Y_\beta^{ss}=Y_\beta\cap S_\beta. \]
Then one proves
\begin{align} \nonumber
  Y_\beta^{ss}&=S_\beta\cap W_\beta\\ \label{bmb}
        &=\{\mu\in S_\beta\;|\;\beta=\beta_\mu\}, 
\end{align}
and
\begin{equation}\label{sky}
  S_\beta=K\cdot Y_\beta^{ss}.
  \end{equation}
Finally,
\begin{equation}\label{wb-sb}
 W_\beta\setminus\{0\}\subset  S_\beta\cup\bigcup_{||\beta'||>||\beta||}S_{\beta'}. 
\end{equation}

\begin{lem}\label{lem:moment-alpha}
  If $M_\mu\in\Lh$, then
  \[ M_\mu=\sum_{ijk}|\mu_{ij}^k|^2\,\alpha_{ij}^k. \]
   \end{lem}

  \Pf For the canonical basis $e_1,\ldots,e_n$ of $\C^n$ we have
  \[ \langle M_\mu e_k,e_k\rangle = - 2\sum_{ij}|\langle L_{e_k}^\mu e_i,e_j\rangle|^2 + \sum_{ij}|\langle L_{e_i}^\mu e_j,e_k\rangle|^2, \]
  so
  \begin{align*}
    M_\mu &=\sum_k\langle M_\mu e_k,e_k\rangle E_{kk} \\
          &=-\sum_{ijk}|\langle L_{e_i}^\mu e_j,e_k\rangle|^2E_{ii}
            -\sum_{ijk}|\langle L_{e_i}^\mu e_j,e_k\rangle|^2E_{jj}
            +\sum_{ijk}|\langle L_{e_i}^\mu e_j,e_k\rangle|^2E_{kk}\\
          &= \sum_{ijk}|\mu_{ij}^k|^2\,\alpha_{ij}^k,
  \end{align*}
  as desired. \EPf

  \begin{cor}\label{ineq}
    If $m(\mu)\in\Lh$ then $m(\mu)\in\mathrm{Conv}\left(\{\alpha_{ij}^k:\mu_{ij}^k\neq0\}\right)$ (convex hull); in particular,
    \[ E(\mu)=||m(\mu)||^2\geq||\beta_\mu||^2, \]
and equality holds if and only if $m(\mu)=\beta_\mu$ if and only if $\mu$ is a soliton,
and in this case $\mu\in S_\beta$ for $\beta\in\mathcal B$ the unique
element of $\bar\Lh^+$ $\mathrm{Ad}_K$-conjugate to $\beta_\mu$.  
    In general, since $E$ is $K$-invariant, from equations~(\ref{bmb})
    and~(\ref{sky})
    we get that
    \begin{equation}\label{eb}
 E(\mu)\geq||\beta||^2\ \mbox{for all $\mu\in S_\beta$,} 
\end{equation}
and equality holds if and only if $\mu$ is a soliton, in which 
case $m(\mu)$ is $\Ad_K$-conjugate to $\beta$ and to $\beta_\mu$. 
  \end{cor}

\begin{cor}
  If $\mu$ is a soliton and $0$ is an eigenvalue of
  $D_\mu$, then $||\beta_\mu||\leq1$.
\end{cor}

\Pf By replacing $\mu$ by a $\mathrm{Ad}_K$-conjugate, we may assume
$\beta_\mu=m(\mu)$, so 
$\beta_\mu+||\beta_\mu||^2I=\frac1{||M_\mu||^2}D_\mu$ has $0$ as
an eigenvalue, that is, $(\beta_\mu,E_{ii})=-||\beta_\mu||^2$
for some $i=1,\ldots,n$. 
Finally, $||\beta_\mu||\geq|(\beta_\mu,E_{ii})|=||\beta_\mu||^2$. \EPf

\medskip

Since the stratification~(\ref{strat})
is $G$-invariant, it naturally induces a stratification of any
$G$-invariant subvariety of $V_n\setminus\{0\}$. 

\begin{prop}\label{stratum-1/n}
The stratum $S_\beta\cap \mathcal J_n$ for $\beta=-\frac1nI$ precisely
consists of the $n$-dimensional semisimple Jordan algebras. 
\end{prop}

\Pf It follows from Proposition~\ref{min-fn} that if 
$\mu\in\mathcal J_n$ is semisimple then there is $g\in G$ 
with $M_{g\cdot\mu}$ a scalar matrix. Therefore
$\mu\in S_\beta\cap\mathcal J_n$. 
Conversely, assume that $\mu\in S_\beta\cap \mathcal J_n$. 
Then the integral curve $\{\mu(t)\}$  of the $(-\nabla E_n)$-flow 
with $\mu(0)=\mu$ is contained in $G\cdot\mu$ and 
converges to a soliton  in $S_\beta$. 
In particular all the eigenvalues of $M_{\mu(t)}$ 
are negative for sufficiently large $t$. In the second half of the 
proof of Proposition~\ref{min-fn}, the argument only needs
this information to imply that $\mu(t)$ must be semisimple for 
sufficiently large $t$. Hence $\mu$ is semisimple, too. \EPf

\subsection{Proofs of Proposition~\ref{deg-decr-e} 
and Theorem~\ref{ss}}

If $\mu\to\nu$ then $G\cdot\nu\subset\overline{G\cdot\mu}$. 
Say $\mu\in S_\beta$ for some $\beta\in\mathcal B$. Then~(\ref{sb})
implies that $\nu\in S_\beta\cup \bigcup_{||\beta'||>||\beta|}S_{\beta'}$.
Owing to Corollary~\ref{ineq}, we have
\[ E(\nu)\geq||\beta||^2=E(\mu), \]
proving Proposition~\ref{deg-decr-e}. 

We move to the proof of the theorem. 
Recall that $S_\beta\cap\mathcal J_n$ for $\beta=-\frac1nI$
precisely consists of the $n$-dimensional semisimple Jordan algebras
(Proposition~\ref{stratum-1/n}), so it is open in
$\mathcal J_n$, say thanks to Albert's criterion for semisimplicity
in terms of the trace form (subsection~\ref{sec:ss}). 
Suppose $\mu$ is semisimple and $\nu\to\mu$.
Then $\nu$ is semisimple. Moreover $\mu\in\overline{G\cdot\nu}$ implies
that $\dim G\cdot\mu\leq\dim G\cdot\nu$ and therefore
$\dim\mathrm{Der}(\mu)\geq\dim\mathrm{Der}(\nu)$. 
We will prove the reverse inequality by 
using a result of Jacobson~\cite{jacobson1949}
asserting that every derivation $D$ of a semisimple Jordan algebra 
$\Ja$ (over a field of characteristic zero) is inner, in the sense 
that it is given as $D=\sum_i[L_{a_i},L_{b_i}]$ for some 
$a_i$, $b_i\in\Ja$. It implies that 
\begin{equation}\label{semi-cont}
  \dim\mathrm{Der}(\mu)=\dim\mathrm{span}\{[L_{x_i}^\mu,L_{x_j}^\mu]\}_{ij}
  \leq\dim\mathrm{span}\{[L_{x_i}^\nu,L_{x_j}^\nu]\}_{ij}=
  \dim\mathrm{Der}(\nu)
\end{equation} 
for a fixed basis $x_1,\ldots,x_n$ of $\C^n$, by lower semicontinuity
of the dimension of the span in terms of the multiplication in $V_n$.

Now $G\cdot\mu$ and $G\cdot\nu$ have the same dimension, which implies that 
$G\cdot\mu=G\cdot\nu$, that is, $\mu$ and $\nu$ are isomorphic. 
This means that $G\cdot\mu$ is closed in $S_\beta\cap\mathcal J_n$
for $\beta=-\frac 1nI$. 
Since there are only finitely many $G$-orbits in
$S_\beta\cap\mathcal J_n$
(Proposition~\ref{stratum-1/n}), they are all open in 
$S_\beta\cap\mathcal J_n$, and hence in $\mathcal J_n$. 
In particular $G\cdot\mu$ is open in $\mathcal J_n$ and 
hence $\mu$ is rigid. This finishes the proof of Theorem~\ref{ss}. 

\section{Low dimensional Jordan algebras and other examples}

The proof of Theorem~\ref{low} regarding Jordan algebras of dimension at most~$4$ is given in Tables~\ref{J1}-\ref{J4bis}. 
In this section we explain how to read them and explain certain cases in more details. 
We start with some remarks of a general nature. 

\subsection{Decomposable algebras}

A Jordan algebra is called \emph{undecomposable} if it is 
not isomorphic to the direct product of two Jordan algebras,
and \emph{decomposable} otherwise.  
The following lemma, whose proof is easy,
shows that we can restrict our search for
solitons to the undecomposable Jordan algebras. 

\begin{lem}
 if $\mu\in V_n$ and $\nu\in V_m$ are solitons 
and $M_\mu=c_\mu I+D_\mu$, $M_\nu=c_\nu I+D_\nu$, 
then $\mu\times c\nu\in V_{n+m}$ is a soliton,
where $c=\sqrt{\frac{c_\mu}{c_\nu}}$.
\end{lem}

\subsection{Unitalization of Jordan algebras}
The following result is an easy check using, say,
Lemma~\ref{lem:moment-alpha}.

\begin{lem}
If a Jordan algebra $\Ja=(\C^n,\mu)$ does not carry a
unit element and we adjoin a unit element to $\Ja$ 
to obtain $\hat\Ja=(\C^{n+1},\hat\mu)$, then the moment matrix
\[ M_{\hat\mu}=\begin{pmatrix}M_\mu& \\ & -(2n+1)\end{pmatrix}. \]
In particular, if $\mu\in\mathcal J_n$ is a soliton
and $M_\mu=c_\mu I+D_\mu$, then 
$\widehat{\sqrt c\mu}\in\mathcal J_{n+1}$ is a soliton, where
$c=\frac{2n+1}{-c_\mu}$, and
\[ M_{\widehat{\sqrt c\mu}}=-(2n+1)I+\begin{pmatrix}c D_\mu& \\ & 0\end{pmatrix}; \]
in this case $E_{n+1}(\widehat{\sqrt c\mu})=\frac{E_n(\mu)}{E_n(\mu)+1}$. 
 \end{lem}

Basic examples are $\widehat{\Ja_{2,2}}=\Ja_{3,4}$,
$\widehat{\Ja_{2,3}}=\Ja_{3,7}$, $\widehat{\Ja_{2,5}}=\Ja_{3,3}$,
$\widehat{\Ja_{3,13}}=\Ja_{4,25}$,
$\widehat{\Ja_{3,17}}=\Ja_{4,39}$ (see the tables in section~\ref{intro}).

\subsection{The regular representation of a semisimple Jordan algebra}

Let $\Js=(\C^n,\mu)$ be a semisimple Jordan algebra, 
and let $\Jn$ be the underlying 
vector space of $\Js$, which we consider as a $\Js$-module under the regular
representation. We put $\Jn^2=0$ so that $\Ja=\Js+\Jn=(\C^{2n},\tilde\mu)$ 
is a Jordan algebra. Choose the Hermitean product such that $\Js\perp\Jn$, 
on $\Js$ it is the Hermitean product that makes $\mu$ a soliton
with $M_\mu=-\frac 1n I$,
and on $\Jn$ it is isometric to $\Js$. Then a simple
calculation shows that $M_{\tilde\mu}=-\frac 1n I=-\frac 2{2n} I$. 
Basic examples are $\widetilde{\Ja_{1,1}}=\Ja_{2,1}$ and 
$\widetilde{\Ja_{2,4}}=\Ja_{4,22}$ (see the tables in section~\ref{intro}). 

\subsection{Jordan algebras in dimensions $1$,  $2$ and~$3$}
Excluding the trivial algebras, there is 
one isomorphism class of complex $1$-dimensional 
Jordan algebras, five isomorphism
classes of complex $2$-dimensional Jordan
algebras, and $19$ isomorphism
classes of complex $3$-dimensional Jordan 
algebras~\cite{kashuba-martin2014};
those are listed in Tables~\ref{J1}, \ref{J2} and~\ref{J3} in soliton form.
One computes the moment matrices using 
Proposition~\ref{moment-jordan}, and then 
uses Proposition~\ref{prop:struct}
to find a soliton in each isomorphism class. It is essentially 
enough to work with non-semisimple and undecomposable algebras. 
The stratification and invariants are collected in 
Tables~\ref{J1bis}, \ref{J2bis} and~\ref{J3bis}.
To give one example in dimension~$3$, 
consider $\Ja_{3,4}$. One replaces the standard basis $e_1$, $e_2$, $n_1$
given in~\cite{kashuba-martin2014} (see Table~2 therein, where the algebra is 
listed as~$\mathcal T_{10}$) 
by $f_1$, $f_2$, $n_1$, where $f_1=e_1+e_2$, $f_2=e_1-e_2$, in order
to diagonalize the moment matrix, and then replaces $f_2$ by 
$f_2'=k f_2$, where $k^4=5/3$, in order to find a soliton.

\subsection{Jordan algebras in dimension~$4$}
Excluding the trivial algebra, there are~$72$
isomorphism classes of complex $4$-dimensional 
Jordan algebras~\cite{kashuba-martin2014}. 
Those are listed in Tables~\ref{J4i}
and \ref{J4ii} in soliton form, except for~$\Ja_{4,63}$. 
The stratification and invariants are collected in 
Table~\ref{J4bis}.
We omit the tedious calculations and only give a few
typical examples.

\subsubsection{The orbit of $\Ja_{4,66}$ is distinguished.}
We search for solitons among isomorphic algebras
of the form $n_1^2=\alpha n_2$, $n_4^2=\beta n_3$, $n_1n_2=\gamma n_3$.
Using Propositions~\ref{moment-jordan} and~\ref{prop:struct}
we find one with $\alpha^2=\beta^2=4$ and $\gamma^2=3$. 

\subsubsection{The orbits of $\Ja_{4,16}$, $\Ja_{4,17}$
  and $\Ja_{4,25}$ are distinguished.}
We compute that the energy of~$\Ja_{4,16}$ 
in the basis given in~\cite[Table 3]{kashuba-martin2014}
is~$27/49$.  
Consider a point of minimum $\nu$
of $E_4$ in~$\overline{ G\cdot\mu}$.
According to~\cite{martin},
the first level degenerations of $\mu$ are $\Ja_{4,28}$,
$\Ja_{4,31}$ and $\Ja_{4,50}$, whose orbits contain solitons
with energy levels respectively $1$, $1$ and $3/4$. Since
$E(\nu)<27/49$, $\nu$ cannot lie in those orbits. Still
according to~\cite{martin}, the only other possible degenerations
of $\mu$ are to~$\Ja_{4,i}$, where $i=48$, $64$, $66$, $67$, $68$, $70$. 
However $\Ja_{4,48}$ has a soliton with energy level $5/6$
and the other ones have solitons with energy level well above $1$,
so again $\nu$ cannot lie in those orbits. The only remaining 
possibility is that $\nu\in G\cdot\mu$. Hence $\Ja_{4,16}$ 
is distinguished. It follows from Corollary~\ref{ineq} that
$E(\nu)=||\beta_\nu||^2=||\beta_\mu||^2=1/2$. 

We compute that the energy of~$\Ja_{4,17}$
and $\Ja_{4,25}$ are~$3/5$ and $5/9$, respectively, and proceed
similarly in those cases. 

In Table~\ref{J4i} we have written the approximate values of the 
structural constants for these solitons, which were obtained 
by using computer software. 

\subsubsection{The orbit of $\Ja_{4,63}$ is not distinguished.}
Denote $(\C^4,\mu)=\mathfrak A_{4,{63}}$
and $(\C^4,\nu)=\mathfrak A_{4,{64}}$. It is known that $\mu\to\nu$. 
Indeed, let $g_t=\left(\begin{smallmatrix} 1&&&\\&t&&\\&&t&\\&&&t\end{smallmatrix}\right)\in G$
and put $\mu_t=g_t^{-1}\cdot\mu$. Then $\lim_{t\to0}\mu_t=\nu$.
Note that $\nu$ is a soliton and put $\beta=m(\nu)=\left(\begin{smallmatrix}-1&&&\\&-\tfrac12&&\\&&0&\\
  &&&\tfrac12\end{smallmatrix}\right)\in\bar\Lh^+$. It is immediate to see that $\mu\in W_\beta$. Since $\inf E_4(G\cdot\mu)=||\beta||^2=\frac32$,
we deduce from~(\ref{sb}), (\ref{wb-sb}) and~(\ref{eb}) that $\mu\in S_\beta$.

Now suppose, to the contrary, that $\lambda\in G\cdot\mu$
is a soliton. Then $m(\lambda)$ is
$\Ad_K$-conjugate to $\beta$ (Corollary~\ref{ineq}),
and by replacing $\lambda$ by an element in its $K$-orbit 
we may assume $m(\lambda)=\beta$. 
 Since $\frac1{||\lambda||^2}D_\lambda=\beta+||\beta||^2I
=\frac12\left(\begin{smallmatrix}1&&&\\&2&&\\&&3&\\&&&4\end{smallmatrix}\right)$,
 this implies that this matrix
is a derivation of~$\lambda$. This matrix has pairwise different eigenvalues,
so $\lambda$ must be given by
\begin{equation}\label{lambda}
 \lambda(n_1,n_1)=a n_2,\ \lambda(n_1,n_2)=b n_3,\ \lambda(n_1,n_3)=c n_4,\ \lambda(n_2,n_2)=d n_4, 
\end{equation}
and the other products zero,  
for some complex constants $a$, $b$, $c$, $d$. 

Suppose now $\lambda$ is given by~(\ref{lambda}). 
We finish by checking
that: (i) $\lambda$ can only be isomorphic to $\mu$ if $a=0$ and 
$b$, $c$, $d\neq0$; (ii) if $a=0$, $\lambda$ can be a soliton only 
if $b=\pm c$ and $d=0$. This will prove that there are no 
solitons in $G\cdot\mu$. 

Write $\Ja=\Ja_{4,63}=(\C^4,\mu)$ and $\Jb=(\C^4,\lambda)$. Note that 
$\Ja^2=\mathrm{span}(n_3,n_4)$ and 
$\Ja^3=\mathrm{span}(n_4)$, so $\Ja/\Ja^3=\mathfrak A_{3,{18}}$.
On the other hand,
$\Jb^2=\mathrm{span}(an_2, bn_3,cn_4, d n_4)$ and 
$\Jb^3=\mathrm{span}(abn_3, bcn_4, adn_4)$. Suppose
$\Ja$ and $\Jb$ are isomorphic. If $a\neq0$
then, owing to $\dim\Jb^2=2$, we have $b=0$ or $c=d=0$. 
In both cases we get $\Jb/\Jb^3=\Ja_{3,19}$. This shows that
$a=0$. Now $\dim\Jb^2=2$ and $\dim\Jb^3=1$ imply that $b\neq0$ and 
$c\neq0$. We also have $d\neq0$, for otherwise $\Jb$ would be 
isomorphic to $\Ja_{4,64}$. This proves~(i).

We turn to (ii). Suppose $a=0$ and $\lambda$ is a soliton.
We compute that 
\begin{equation}\label{ml} 
M_\lambda=\begin{pmatrix}-2b^2-2c^2&&&\\
&-2b^2-2d^2&&\\
&&-2c^2+2b^2&\\
&&&2c^2+d^2\end{pmatrix}. 
\end{equation}
Since $M_\lambda$ and $\beta$ are multiples, this 
immediately gives that $b^2=c^2$ and $d=0$,
which proves~(ii).

\section{Closed $SL(n,\C)$-orbits in $\mathcal J_n$}

Since there are no closed $GL(n,\C)$-orbits in $\mathcal J_n$ (not even in 
$V_n\setminus\{0\}$), precisely because of multiples of the identity,
it arises the natural question of knowing whether 
the subgroup $SL(n,\C)$ admits closed orbits in $\mathcal J_n$. 

\begin{prop}
Let $\mu\in\mathcal J_n$, $\mu\neq0$. Then:
\begin{enumerate}
\item The orbit $SL(n,\C)\cdot\mu$ is closed if and only if $\mu$ 
is a semisimple Jordan algebra. 
\item If $\mu$ is not semisimple then $0$ lies in the closure 
of  $SL(n,\C)\cdot\mu$.
\end{enumerate}
\end{prop}

\Pf The orbit $SL(n,\C)\cdot\mu$ is closed if and only the moment map 
$m^{SL_n}$ of the $SL(n,\C)$-action on $V_n$ vanishes at some point 
$\nu\in SL(n,\C)\cdot\mu$. 
Since $m^{SL_n}$ is obtained from $m^{GL_n}$ by pos-composing 
with the projection $i\mathfrak{u}(n)=\R\oplus i\mathfrak{su}(n)
\to i\mathfrak{su}(n)$, the latter is equivalent to $m^{GL_n}(\nu)=-\frac1nI$,
which means that $\nu$ and $\mu$ are semisimple, owing to Proposition~\ref{min-fn}.

Suppose now $\mu$ is not semisimple. Then $SL(n,\C)\cdot\mu$ 
is not closed
and $\overline{SL(n,\C)\cdot\mu}$ contains a closed orbit, say, 
$SL(n,\C)\cdot\nu$. If $\nu\neq0$ then $\nu$ is semisimple by 
part~(a). 
But then $\mu\to\nu$ and this contradicts the rigidity of semisimple 
Jordan algebras. \EPf

\section{Partial results, open problems and conjectures}

It is interesting to note that the application of GIT to the study of 
(commutative) Jordan algebras has many similarities with the case of
(anti-commutative) Lie algebras. However the Jordan identity (in each of 
its disguises) seems to be
more difficult to use than the Jacobi identity. 
In particular, for Jordan algebras in general the left multiplications are not
derivations of the algebra, in flagrant contrast with Lie algebras.
So proofs of results for Lie algebras which depend on this property cannot be
simply carried over to the context of Jordan algebras, and throughout this work 
we have tried to find alternative lines of arguments, with some success.
The partial results that we collect in this section are somehow related
to this situation.

Unless explicitly stated, throughout this section we 
let $\mu\in\mathcal J_n$ be a Jordan soliton,
$\Ja=(\C^n,\mu)$. We write  the moment matrix $M_\mu=c_\mu I+D_\mu$,
where $D_\mu$ is a Hermitean derivation, according to 
Proposition~\ref{prop:struct}. Let also $\Jn$ denote the 
radical of $\Ja$. 

\subsection{The annihilator}\label{ann}
The annihilator of $\Ja$ is $\mathrm{Ann}(\Ja)=\{x\in\Ja:L^\mu_x=0\}$.
It is clear from $L^\mu_{D_\mu x}=[D_\mu,L^\mu_x]$ that $D_\mu$ preserves 
$\mathrm{Ann}(\Ja)$, and it follows from
Corollary~\ref{cor:moment-jordan} that the eigenvalues of $D_\mu$ on 
$\mathrm{Ann}(\Ja)$ are positive.

\subsection{Basic calculation}
Let $x$ be an eigenvector of $D_\mu$ with eigenvalue $d$. 
Then
\begin{align} \nonumber
(M_\mu,[L^\mu_x,L^{\mu*}_x]) & = \mathrm{Tr}(M_\mu[L^\mu_x,L^{\mu*}_x])\\ \nonumber 
& = \mathrm{Tr}(D_\mu[L^\mu_x,L^{\mu*}_x])\\ \label{mm-basic} 
& = \mathrm{Tr}([D_\mu, L^\mu_x]L^{\mu*}_x])\\ \nonumber
& = \mathrm{Tr}(L^\mu_{D_\mu x}L^{\mu*}_x)\\ \nonumber
&= d||L^\mu_x||^2. \nonumber
\end{align}
On the other hand, due to~(\ref{moment-map}) 
the left-hand side of~(\ref{mm-basic}) also equals
$\langle [L^\mu_x,L^{\mu*}_x]\cdot\mu,\mu\rangle$,
so we deduce
\begin{equation}\label{basic} 
d||L^\mu_x||^2 = ||L^{\mu*}_x\cdot\mu||^2-||L^\mu_x\cdot\mu||^2. 
\end{equation}
We use this formula and some variations below.

\begin{lem}\label{lem:ort}
Let $\mu\in\mathcal J_n$ be a soliton, and 
let $x$, $y\in\C^n$ be eigenvectors of $D_\mu$ with 
corresponding (real) eigenvalues $d_x$, $d_y$. 
If $d_x\neq d_y$ then $(L_x^\mu,L_y^\mu)=0$. 
\end{lem}

\Pf The basic calculation yields
\[ d_x(L^\mu_x,L^\mu_y)=(M_\mu,[L^\mu_x,L^{\mu*}_y])=
\langle L^{\mu*}_y\cdot\mu,L^{\mu*}_x\cdot\mu\rangle -\langle L^\mu_x\cdot\mu,L_y^\mu\cdot\mu\rangle. \]
We interchange $x$ and $y$ in these equations to obtain that 
\[ d_x(L^\mu_x,L^\mu_y) = \overline{d_y(L^\mu_y,L^\mu_x)} = d_y(L^\mu_x,L^\mu_y), \]
which proves the desired result. \EPf

\subsection{The kernel of $D_\mu$}
Proposition~\ref{prop:struct}(b) shows that the kernel of~$D_\mu$ 
contains a maximal semisimple subalgebra of $\Ja=(\C^n,\mu)$.

\begin{quest}\label{q1}
For a soliton $\mu\in\mathcal J_n$ 
is it true that $\ker D_\mu$ is a maximal semisimple subalgebra?
In other words, are the eigenvalues of $D_\mu$ restricted to~$\Jn$
different from zero? 
\end{quest}

\subsection{Positivity of $D_\mu$ on $\Jn$}
A positive answer to Question~\ref{q2} implies a positive answer to 
Question~\ref{q1}.

\begin{quest}\label{q2}
For a soliton $\mu\in\mathcal J_n$,
are the eigenvalues of $D_\mu$ on $\Jn$ positive? 
\end{quest}

For soliton Jordan algebras $\Ja$ satisfying $\Ja^3=0$, trivially
all left-multiplications
are derivations of the algebra, and we can use 
a standard argument to answer yes to Question~\ref{q2}.

\begin{prop}
If the soliton $\Ja=(\C^n,\mu)$ satisfies $\Ja^3=0$, then
all eigenvalues of $D_\mu$ are positive.
\end{prop}

\Pf The assumption $\Ja^3=0$ implies that $\Ja$ is a nilpotent 
Jordan algebra and hence $L_x^\mu$ is a nilpotent operator 
for all $x\in\Ja$. 

Let $x\in\C^n$ be an eigenvector of $D_\mu$ with corresponding 
eigenvalue $d\in\R$. Since $L_x^\mu\cdot\mu=0$, the basic 
calculation~(\ref{basic}) 
immediately gives $d\geq0$. If in addition $d=0$, then 
$L^{\mu*}_x$ is also a derivation of $\mu$, thus 
$[L^{\mu*}_x, L^\mu_x]=L^\mu_{L_x^{\mu*}x}$ by the definining 
condition of a derivation. 

Since $D_\mu$ is Hermitean, the orthogononal 
decomposition $\Ja=\Ja^2\oplus\Ja^{2\perp}$ is $D_\mu$-invariant, so we may 
assume $x\in\Ja^2$ (resp.~$x\in\Ja^{2\perp}$). In any case
$\langle L_x^{\mu*}x,y\rangle =\langle x,xy\rangle=0$ for 
all $y\in\Ja$, as $xy\in\Ja^3=0$ (resp.~$xy\in\Ja^2$). This 
shows $L_x^{\mu*}x=0$.

Now $L_x^\mu$ is a normal and nilpotent operator, hence
$x\in\mathrm{Ann}(\Ja)$. The result in subsection~\ref{ann}
contradicts our assumption that $d=0$. Hence $d>0$. \EPf

\medskip

For non-associative Jordan algebras, we have the following 
partial result. 

\begin{prop}
Let $x$, $y\in\C^n$ be eigenvectors of $D_\mu$ with corresponding eigenvalues
$d_x$, $d_y$. If $[L^\mu_x,L_y^\mu]\neq0$ then $d_x+d_y\geq0$. 
\end{prop}

\Pf We consider $D:= [L^\mu_x,L_y^\mu]$ and compute   
  \begin{align*}
    [D_\mu,D]&=[[D_\mu,L^\mu_x],L^\mu_y]+[L^\mu_x,[D_\mu,L^\mu_y]]\\
     & = [L^\mu_{D_\mu x},L^\mu_y]+[L^\mu_x,L^\mu_{D_\mu} y]\\
             &= (d_x+d_y)D.
          \end{align*}
We proceed as in the basic calculation to get
\[ (d_x+d_y)||D||^2=\mathrm{Tr}([D_\mu,D]D^*)=\mathrm{Tr}(D_\mu[D,D^*])=
\mathrm{Tr}(M_\mu[D,D^*])=||D^*\cdot\mu||^2, \]
since $D$ is a (inner) derivation, and the result follows. \EPf

\medskip

If $\Jn=(\C^n,\mu)$ 
is a nilpotent Jordan soliton generated by a single element $n$, 
then $n$ and its powers must be eigenvectors of $D_\mu$.  
Owing to $\mathrm{Tr}(D_\mu)=\frac{\mathrm{Tr} D_\mu^2}{-c_\mu}>0$, 
we also get $D_\mu>0$
in this case.

\subsection{Orthogonality of $\Js$ and $\Jn$}

Since $D_\mu$ is Hermitean, the orthogonal complement of $\Jn$ 
with respect to the Hermitean product is $D_\mu$-invariant. 

\begin{quest}\label{q5}
For a soliton $\mu\in\mathcal J_n$, $\Ja=(\C^n,\mu)$, is $\Jn^\perp$ a 
semisimple subalgebra of $\Ja$? 
\end{quest}

If the answer to Question~\ref{q5} is yes, then $\Jn^\perp$ will be 
maximal semisimple.
A positive answer to Quesion~\ref{q1} implies a positive answer to 
Question~\ref{q5}.

\subsection{Reduction to nilpotent Jordan algebras}

\begin{quest}
  If $\Ja=(\C^n,\mu)$ is a Jordan soliton and $\Jn$ is the radical of $\Ja$,
  is it true that $\mu|_{\mathfrak N\times\mathfrak N}$ is a soliton?
  Conversely, given a nilpotent Jordan soliton $\Jn$ and a semisimple
  Jordan algebra $\Js$ such that  $\Ja=\Js+\Jn$ is a Jordan 
  algebra, can we extend the metric from $\Jn$ to $\Ja$ so that
  $\Ja$ becomes a soliton? 
\end{quest}

 
\providecommand{\bysame}{\leavevmode\hbox to3em{\hrulefill}\thinspace}
\providecommand{\MR}{\relax\ifhmode\unskip\space\fi MR }
\providecommand{\MRhref}[2]{%
  \href{http://www.ams.org/mathscinet-getitem?mr=#1}{#2}
}
\providecommand{\href}[2]{#2}

\end{document}